\documentclass[12pt]{amsart}
\usepackage{amssymb, amsmath, amsfonts}
\usepackage[raiselinks=false,colorlinks=true,
citecolor=blue,urlcolor=blue,linkcolor=blue,
bookmarksopen=true,pdftex]{hyperref}
\newtheorem{theorem}{Theorem}[section]
\newtheorem{proposition}[theorem]{Proposition}
\newtheorem{lemma}[theorem]{Lemma}
\newtheorem{remark}[theorem]{Remark}

\newtheorem{example}[theorem]{Example}

\newcommand{\bc}{\mathbb{C}}

\renewcommand{\hom}{\textrm{Hom}}

\newcommand{\SO}{\textrm{SO}}

\newcommand{\Fix}{\text{Fix}}

\newcommand{\res}{\textrm{Res}}
\newcommand{\diag}{\textrm{diag}}

\begin{document}
\baselineskip=15.5pt
\title[Cohomology classes in uniform lattices of $\SO(n,\mathbb{H})$]{Non-vanishing cohomology 
classes in uniform lattices of $\SO(n,\mathbb{H})$  and automorphic representations}  
\author[A. Mondal]{Arghya Mondal}
\author[P. Sankaran]{Parameswaran Sankaran}
\address{The Institute of Mathematical Sciences, CIT
Campus, Taramani, Chennai 600113, India}
\email{arghya@imsc.res.in}
\email{sankaran@imsc.res.in}
\subjclass[2010]{22E40, 22E46\\
Keywords and phrases: Geometric cycles, automorphic representation.}
\date{}

\begin{abstract}   
Let $X$ denote the non-compact globally Hermitian symmetric space of type $DIII$, namely, $\SO(n,\mathbb{H})/U(n)$.
Let $\Lambda$ be a uniform torsionless lattice in $SO(n,\mathbb{H})$.  
In this note 
we construct certain complex analytic submanifolds in the locally symmetric space $X_\Gamma:=\Gamma\backslash \SO(n,\mathbb{H})/U(n)$ for certain finite index sub lattices $\Gamma\subset \Lambda$ and show that their dual cohomology classes in $H^*(X_\Gamma;\mathbb{C})$ are not in the image of the 
Matsushima homomorphism $H^*(X_u; \mathbb{C})\to H^*(X_\Gamma;\mathbb{C})$, where $X_u=\SO(2n)/U(n)$ is the compact dual of $X$.   These submanifold arise as sub-locally symmetric spaces which are totally geodesic, and, when $\Lambda$ satisfies 
certain additional conditions, they are non-vanishing `special cycles'.  
Using the fact that $X_\Lambda$ is a K\"ahler manifold, we deduce the occurrence in $L^2(\Lambda\backslash SO(n,\mathbb{H})$ of a certain irreducible    
representation $(\mathcal{A}_\mathfrak{q}, A_\mathfrak{q})$ with non-zero multiplicity when $n\ge 9$.  
The representation $\mathcal{A}_\mathfrak{q}$ is associated to a certain $\theta$-stable parabolic subalgebra $\mathfrak{q}$ of $\mathfrak{g}_0:=\mathfrak{so}(n,\mathbb{H})$.   Denoting the smooth $U(n)$-finite vectors of $A_{\mathfrak{q}}$ by $A_{\mathfrak{q},U(n)}$, the representation $\mathcal{A}_\mathfrak{q}$ is characterised by the 
property that $H^{p,p}(\mathfrak{g}_0\otimes\mathbb{C},U(n); A_{\mathfrak{q},U(n)})\cong H^{p-n+2,p-n+2}(\SO(2n-2)/U(n-1);\mathbb{C}),~p\ge 0$, for  $n\ge 9$. 
\end{abstract}
\maketitle
\section{Introduction}\label{intro}
Let $G$ be a connected non-compact real semi simple Lie group with finite centre and let $K\subset G$ be a maximal compact subgroup.  Denote by $\mathfrak{g}_0$ the Lie algebra of $G$ and by $\mathfrak{g}$ its complexification.  The symbols 
$\mathfrak{k}_0, \mathfrak{k}$, likewise denote the Lie algebra of $K$ and its complexification respectively.  Let $\theta:\frak{g}_0\to \frak{g}_0$ denote the Cartan involution that fixes $\frak{k}_0$, so that $\frak{g}_0=\frak{k}_0\oplus \frak{p}_0$ where $\frak{p}_0$ is 
the $-1$-eigenspace of $\theta$.  We denote by $\frak{p}$ the complexification of $\frak{p}_0$ so that 
$\frak{g}=\frak{k}\oplus \frak{p}$.  The Killing form restricted to $\frak{p}_0$ is a positive definite inner product and 
defines a Riemannian metric on $X=G/K$ making it a globally symmetric space. 
If $\Gamma$ is a torsionless  
discrete subgroup of $G$, then the locally symmetric space $X_\Gamma:=\Gamma\backslash X=\Gamma \backslash G/K$ is a 
smooth aspherical manifold.  

Throughout this paper we assume that $\Gamma$ is a uniform lattice in $G=\SO(n,\mathbb{H})\cong SO^*(2n), n\ge 4$.    (Note that $\Gamma$ is irreducible since $G$ is simple.)  

Our aim is to construct closed oriented totally geodesic submanifolds of $X_\Gamma$, 
which are known as {\it special cycles}, such that their Poincar\'e duals are non-zero cohomology classes in $H^*(X_\Gamma;\mathbb{R})$ in the case when $G=\SO(n,\mathbb{H})$. 

They arise as sub locally symmetric spaces  $C(\sigma)=X(\sigma)_{\Gamma(\sigma)}$ 
where $X(\sigma)=G(\sigma)/K(\sigma)$,  $G(\sigma)\subset G$ is the fixed point of an algebraically defined involutive automorphism  
$\sigma$ of $G$ so that $\Gamma(\sigma)=\Gamma\cap G(\sigma)$ is again a lattice in $G(\sigma)$.    

When $C(\sigma)$ is a complex analytic submanifold of $X_\Gamma$, then, in view of the fact that $X_\Gamma$ is a K\"ahler manifold---it is even a complex projective variety---it follows 
immediately that $[C(\sigma)]$ is non-zero.  
In the more general case when $C(\sigma)$ is not complex analytic,  
to show that the Poincar\'e dual to the fundamental class of $X(\sigma)_{\Gamma(\sigma)}\hookrightarrow X_\Gamma$ is non-zero, 
one constructs another special cycle $C(\tau)$ of complementary dimension and shows that the cup product 
$[C(\sigma)].[C(\tau)]$ in $H^*(X_\Gamma;\mathbb{R})$ is non-zero. The existence of such pairs of 
special cycles of complementary dimensions imply, not only the non-vanishing of their Poincar\'e duals, but also the stronger condition that their Poincar\'e duals are cohomology classes that cannot be  representable by $G$-invariant forms on $G/K$.  Equivalently, these cohomology classes are not in the image of the Matsushima homomorphism \cite{matsushima63} $H^*(X_u;\mathbb{C})\to H^*(X_\Gamma;\mathbb{C})$ where $X_u$ denotes the compact dual of $X$.  
The method of constructing such pairs of special 
cycles are well understood; see \cite{mr} and \cite{rs}.

Special cycles were first constructed by  
Millson and Raghunathan \cite{mr}    
in the case of $G=SU(p,q), SO_0(p,q), Sp(p,q)$.  Schwermer and Waldner \cite{sw} considered the case  $G=SU^*(2n)$, and  Waldner \cite{waldner} the case 
when $G$ is the non-compact real form of the exceptional group $G_2$.  Recently the groups $SL(n,\mathbb{R}), SL(n,\mathbb{C}), n\ge 3,$ have been considered 
by Schimpf \cite{schimpf}.  In all these cases, the choice of the uniform lattice under consideration had to be suitably restricted.

We too need to restrict the class of lattices in $G=\SO(n,\mathbb{H})$.   
Let $F$ be a totally real number field and let $\mathcal{O}_F$ denote its ring of integers. 
Let $\mathbf{G}$ be an $F$-algebraic group such that the $\mathbb{R}$-points $G(\mathbb{R})$ is isomorphic as a Lie group to $G=\SO(n,\mathbb{H})$ and such that the $\mathcal{O}_F$-points $G(\mathcal{O}_F)$ is a {\it uniform}  
lattice in $G(\mathbb{R})$.  (Up to commensurability 
any uniform lattice of $G$ arises in this manner, thanks to the arithmeticity theorem of Margulis.)
We shall refer to such an $F$-structure on $G$ as 
type $DIII_u$.

Our basic hypothesis is a certain condition (see Theorem \ref{f-cartan}) on the $F$-structure which guarantees the existence of a Cartan 
involution $\theta=\theta_\mathbb{R}$ induced by 
an $F$-involution $\theta_F$ on $\mathbf{G}$ given by 
conjugation by a {\it diagonal} matrix in $G(F)$. 
Such an $F$-rational Cartan involution will be said to be of {\it 
diagonal type}.    

Suppose that $\Gamma$ is commensurable with $G(\mathcal{O}_F)$.  
Whenever required, $\Gamma$ may be replaced by a finite index subgroup which is stable under finitely many pairwise commuting collection of involutions $\sigma=\sigma_\mathbb{R}$ induced by 
$F$- involutions $\sigma_F$, including a Cartan involution $\theta$ defined over $F$ of diagonal type.  
Once the former condition is met, the latter is readily achieved---one need only consider the intersection 
of $\Gamma$ with the members of the 
collection $\sigma(\Gamma)$ as $\sigma$ varies in the finite collection of commuting involutions.

Our first main result is the following:

\begin{theorem} \label{main1}    Let $n\ge 8$.
Suppose that $\theta$ is an $F$-rational Cartan involution on $G=\SO(n,\mathbb{H})$ of diagonal type where the $F$-structure on $G$ is of type $DIII_u$. 
Suppose that $\Lambda$ is any torsionless lattice in  $\SO(n,\mathbb{H}), n\ge 4,$ that is commensurable with $G(\mathcal{O}_F).$  Then there exist cohomology 
classes in $H^{2p}(X_\Lambda)$ which are not in the 
image of the Matsushima homomorphism $H^*(X_u;\mathbb{C})\to H^*(X_\Lambda)$, in all even dimensions 
$2p$ between $2n-4$ and $n(n-1)-(2n-4)$.  Moreover, 
when $n\equiv 1 ~\mod 2$, there exists a cofinal family $\{\Gamma\}$ of finite index 
subgroups $\Gamma \hookrightarrow\Lambda$ such that 
$H^{n(n-1)/2}(X_\Gamma;\mathbb{C})$ contains classes 
not belonging to the image of the Matsushima homomorphism.  
\end{theorem}

In order to state our next main result, which may be viewed as an application  
of the above theorem, we need to recall 
certain facts about $L^2(\Gamma\backslash G)$ as well as relative Lie algebra cohomology of 
unitary representations. (Here $\Gamma$ is a uniform lattice.)  
To a Haar measure on $G$ is associated in a natural 
manner, a $G$-invariant measure on $\Gamma \backslash G$.  The 
Hilbert space  $L^2(\Gamma\backslash G)$ affords a unitary representation of $G$ (via the right action of $G$ on 
$\Gamma \backslash G$).  From the work of Gelfand and 
Pyatetskii-Shapiro (\cite[Ch. 1, \S2.3]{ggp}, \cite{gp}), $L^2(\Gamma \backslash G)$ is a completed direct sum of irreducible unitary $G$-representations $(\pi, H_\pi)$ each of which occurs with finite multiplicity $m(\pi, \Gamma)$.  The representations of $\pi$ with $m(\pi,\Gamma)>0$ 
are called {\it automorphic} representations. 
Denoting the space of smooth $K$-finite vectors of a $G$-representation $V$ by $V_K$, one has 
the Matsushima isomorphism $H^*(X_\Gamma;\mathbb{C}) 
=H^*(\mathfrak{g},K; L^2(\Gamma\backslash G)_K)
\cong \oplus_{\pi}m(\pi,\Gamma) H^*(\mathfrak{g},K;H_{\pi,K})$.  See \cite{matsushima}.
By a result of D. Wigner $H^*(\mathfrak{g},K; H_{\pi,K})$ is 
non-zero only if, the infinitesimal character of $\pi$ equals that of the trivial representation of $G$.  The representations $(\pi, H_\pi)$ 
with non-vanishing $(\mathfrak{g},K)$-cohomology 
have been classified (up to unitary equivalence) in terms of the $\theta$-stable parabolic subalgebras $\mathfrak{q}$ of $\mathfrak{g}$.  Denote 
the representation corresponding to $\mathfrak{q}$ by $(\mathcal{A}_{\mathfrak{q}},A_\mathfrak{q})$.  If $H^*(\mathfrak{g},K; H_{\pi,K})\ne 0$ then $(\pi,H_\pi)$ is unitarily equivalent to $(\mathcal{A}_\mathfrak{q},A_\mathfrak{q})$ for some $\theta$-stable 
parabolic subgroup of $\mathfrak{g}_0$.

When $X$ is a Hermitian symmetric space, $X_\Gamma$ is a K\"ahler manifold, which is compact since 
$\Gamma$
is uniform.  One has the Hodge decomposition $H^r(X_\Gamma; \mathbb{C})=\oplus_{p+q=r}H^{p,q}(X_\Gamma;\mathbb{C})$.  Also one has the Hodge decomposition $H^r(\mathfrak{g},K;V_K)=\oplus_{p+q=r}H^{p,q} (\mathfrak{g},K;V_K)$ where $V$ is any unitary $G$-representation; see \cite[Ch. II \S4]{bw}.      
Corresponding to each $\mathfrak{q}$, there exists a pair of numbers, $(R_+(\mathfrak{q}),R_-(\mathfrak{q}))$, called the type of $\mathfrak{q}$, 
such that $H^{p,q}(\mathfrak{g}, K;A_{\mathfrak{q},K})=0$ unless $p\ge R_+(\mathfrak{q}), q\ge R_-(\mathfrak{q}), p-q=R_+(\mathfrak{q})-R_-(\mathfrak{q}).$ Moreover  $H^{R_+(\mathfrak{q}),R_-(\mathfrak{q})}(\mathfrak{g},K;A_{\mathfrak{q},K})\cong \mathbb{C}$.  

We shall show, when $G,\Gamma$ are as in Theorem \ref{main1},  that, for $n\ge 8$ there exists a unitary representation $\mathcal{A}_\mathfrak{q}$, where the 
$\theta$-stable parabolic subalgebra 
$\mathfrak{q}$ is of type $(n-2,n-2)$, and {\it none} of type $(n-l,n-l)$ for $2<l<n$ or $l=1$.  This representation is {\it unique} up to 
unitary equivalence when $n\ge 9$. 

It is an important open problem in the study of cohomology of locally symmetric spaces, to determine, for a 
given $\Gamma\subset G$ and $\mathfrak{q}$, whether $m(\mathcal{A}_\mathfrak{q},\Gamma)$ is positive, let alone 
its value.  We have the following result.  

\begin{theorem} \label{main2}   Let $n>8$.  Let $\mathbf{G}$ be an $F$-structure on $\SO(n,\mathbb{H})$ of type 
$DIII_u$ admitting an $F$-rational Cartan involution as in Theorem \ref{main1}.  
Then, for any torsionless lattice $\Gamma$ commensurable with $G(\mathcal{O}_F)$ we have 
$m(\mathcal{A}_\mathfrak{q},\Gamma)> 0$ where 
$\mathfrak{q}$ is of type $(n-2,n-2)$.  
\end{theorem}

The above result seems to be a new addition to the vast body of work on the occurrence of $H_\pi$ in $L^2(\Gamma\backslash G)$ as well as the asymptotics of $m(\pi,\Gamma)$ in various settings (not restricted to the case of $SO(n,\mathbb{H})$) such as when 
$\pi=\mathcal{A}_\mathfrak{q}$ is holomorphic (i.e. of type $(p,0)$) due to Anderson \cite{anderson}, or a discrete series representation due to Clozel \cite{clozel}, DeGeorge and Wallach \cite{degeorge-wallach}, 
or when $\pi=\mathcal{A}_\mathfrak{q}$ where the real reductive subgroup $L\subset G$ with $\mathfrak{l}= 
\mathfrak{q}\cap \overline{\mathfrak{q}}$ is isomorphic to a group of the form 
$(\mathbb{S}^1)^r\times SO(n-r,\mathbb{H})$ due to Li \cite{li}.  See also \cite[Chapter VIII]{bw}, \cite[\S6]{parthasarathy80} and \cite{kumaresan}.
 
The paper is organised as follows.   In \S \ref{arthmeticlattices} we consider $F$-structures on $\SO(n,\mathbb{H})$, construct 
$F$-involutions, including $F$-rational Cartan involution, under certain hypotheses on the $F$-structure, and 
construct special cycles, in \S \ref{specialcycle}, associated to these involutions.
In \S \ref{thetastable} we consider $\theta$-stable parabolic subalgebras of $\SO(n,\mathbb{H})$ and compute their relative Lie algebra cohomology in terms of the certain combinatorial data associated to them.  We also classify, in Proposition \ref{minimalparabolics}, $\theta$-stable parabolic subalgebras of $\mathfrak{so}(n,\mathbb{H})$ of type $(p,q)$ when $p,q \le n-1$.  
Theorems \ref{main1} and \ref{main2} will be proved in \S \ref{proofs}.

\section{Arithmetic lattices in $\SO(n,\mathbb{H})$}\label{arthmeticlattices}
\subsection{$\SO(n,\mathbb{H})$ and $\SO(2n,\mathbb{C})$} \label{realform}
Let $\mathbb{H}$ denote the division algebra of real quaternions with standard $\mathbb{R}$-basis $1, i,j, k=ij$.   We denote by $\tau_c: \mathbb{H}\to \mathbb{H}$ 
the canonical conjugation defined as $\tau_c(q)=\bar{q}:=q_0-q_1i-q_2j-q_3k $ where $q=q_0+q_1i+q_2j+q_3k, q_r\in \mathbb{R}$.  We denote by $N(q)$ the norm of $q$, namely, $N(q)=q.\bar{q}=q_0^2+q_1^2+q_2^2+q_3^2$.  The conjugation $\tau_c$ is an anti-automorphism of the division algebra $\mathbb{H}$. 
One has also the anti-automorphism $\tau_r:\mathbb{H}\to \mathbb{H}$ defined as $\tau_r(q)=q_0+q_1i-q_2j+q_3k=j\tau_c(q)j^{-1}$, known as the {\it reversion}.  Applying $\tau_r$ to each entry we obtain an $\mathbb{R}$-vector space 
automorphism $M_n(\mathbb{H})\to M_n(\mathbb{H})$ again denoted $\tau_r$.  Explicitly, $\tau_r(Q)=Q_0+Q_1i-Q_2j+Q_3k$ where $Q=Q_0+Q_1i+Q_2j+Q_3k, Q_i\in M_n(\mathbb{R})$.  Since transposition $P\mapsto{}^tP$ is an $\mathbb{R}$-algebra 
anti-automorphism of $M_n(\mathbb{R})$, and since $\tau_r$ is an anti-automorphism of $\mathbb{H}$, it is follows that $Q\mapsto {}^t\tau_r(Q)$ is an $\mathbb{R}$-algebra 
anti-automorphism of $M_n(\mathbb{H})$.  

Note that $\mathbb{H}$ splits over $\mathbb{C}$, that is, we have an isomorphism of $\mathbb{C}$-algebras $\mathbb{H}_\mathbb{C}:=\mathbb{H}\otimes_\mathbb{R}\mathbb{C}\cong M_2(\mathbb{C})$, where $ q\otimes 1=q=z+wj\in 
\mathbb{H}, z, w\in \mathbb{C}$ is sent to $\left(\begin{smallmatrix}z & w\\ -\bar{w} & \bar{z}\end{smallmatrix}\right) $.  
Under this isomorphism, $\tau_r(q)$ maps to $\left(\begin{smallmatrix}z &-\bar{w}\\w & \bar{z}\end{smallmatrix}\right)$.  

The inclusion  $\mathbb{H}\hookrightarrow \mathbb{H}_\bc\cong M_2(\mathbb{C})$ yields an obvious embedding 
$M_n(\mathbb{H})\hookrightarrow M_{2n}(\mathbb{C})$ where each quaternion in the domain is viewed as 
a $2\times 2$ matrix with entries in $\mathbb{C}$.   However, it is more convenient for our purposes to use 
a different embedding, namely, $\psi: M_n(\mathbb{H})\to M_{2n}(\mathbb{C})$ defined as $Q\mapsto \left ( \begin{smallmatrix}Z & W \\-\bar{W} & \bar{Z}\end{smallmatrix}\right) $ where $Z, W\in M_{n},(\mathbb{C})$ are defined by $Q=Z+Wj$.   It is easy to see that $\psi$ is an $\mathbb{R}$-algebra homomorphism.   
Note that $\psi({}^tQ)\ne {}^t\psi(Q)$, indeed it is readily verified that $\psi({}^t\tau_r(Q))={}^t\psi(Q)$. 
 This shows that the anti-automorphism ${}^t\tau_r$ on $M_n(\mathbb{H})$ corresponds, under $\psi$, to transposition 
 in $M_{2n}(\mathbb{C})$. 
In particular ${}^t\tau_r(Q)=Q$ if and only if 
$\psi(Q)$ is symmetric if and only if $Z$ is symmetric and $W$ is skew-hermitian.  Similarly ${}^t\tau_r(Q)=-Q$ if 
and only if $Z$ is skew-symmetric and $W$, hermitian. 

The group  $\SO(n, \mathbb{H})$ is defined as $\{Q\in \textrm{SL}(n,\mathbb{H})\mid  {}^t\tau_r(Q)Q=I_n\}$ where $I_n$ 
denotes the $n\times n$ identity matrix.  
By the above discussion, it is clear that $\psi$ restricts to a monomorphism of Lie groups $\SO(n,\mathbb{H})\to \SO(2n,\mathbb{C}),$ again denoted $\psi$.  We shall identify $\SO(n,\mathbb{H})$ with its image under $\psi$.  
The Lie algebra of $\SO(n,\mathbb{H})$ is seen to be $\mathfrak{so}(n,\mathbb{H})=
\{Q\in M_n(\mathbb{H})\mid {}^t\tau_r(Q)+Q=0\}=\{ \left(\begin{smallmatrix}Z & W \\-\bar{W} & \bar{Z}\end{smallmatrix}
\right) \in M_{2n}(\mathbb{C})\mid Z+{}^tZ=0, W-{}^t\bar{W}=0\}$.  It is clear from this description that $\mathfrak{so}(n,\mathbb{H})\otimes_\mathbb{R}\mathbb{C}
=\frak{so}(2n,\mathbb{C}),$ the Lie algebra of $\SO(n,\mathbb{C})$.  
Thus $\SO(n,\mathbb{H})$ is a real form of 
$\SO(2n,\mathbb{C})$. The subgroup $K=\{Z+Wj\in \SO(n,\mathbb{H})\mid Z, W\in M_n(\mathbb{R})\}$ is isomorphic to 
$U(n)$ since ${}^t\tau_r(Z+Wj)={}^tZ-{}^tWj={}^t\overline{(Z+Wj)}$ if $Z,W\in M_n(\mathbb{R})$ where bar denotes 
conjugation in $M_n(\mathbb{R}+\mathbb{R}j) \cong M_n(\mathbb{C})$.   We also have 
$\mathfrak{k}_0=\{A+Bj\in \mathfrak{so}(n,\mathbb{H})\mid A,B\in M_n(\mathbb{R})\}$, which corresponds 
under $\psi$ to 
$\{\left(\begin{smallmatrix}A & B \\-B& A\end{smallmatrix}\right)\mid 
A, B\in M_n(\mathbb{R}), A+{}^tA=0, B={}^tB\}$.   The subgroup $K$ is a maximal compact subgroup of 
$\SO(n,\mathbb{H})$ and will be referred to as the {\it standard} maximal compact subgroup of $\SO(n,\mathbb{H})$.  The 
corresponding Cartan involution $\theta: \SO(n,\mathbb{H})\to \SO(n,\mathbb{H})$ is $Q\mapsto JQJ^{-1}$, 
where $J=jI_n\in M_n(\mathbb{H})$ and its differential at the identity element, again denoted by the same symbol, is  
$\theta:\frak{so}(n,\mathbb{H})\to \frak{so}(n,\mathbb{H})$ is again conjugation by $J$.
Thus $\frak{p}_0$, the $-1$ eigenspace of $\theta$, equals $\{Ai+Bk\in \mathfrak{so}(n,\mathbb{H})\mid  A,B\in M_n(\mathbb{R})\}$ and corresponds 
under $\psi$ to the subspace 
$\{i\left (\begin{smallmatrix}A & B \\B & -A\end{smallmatrix}\right)\in \psi(\mathfrak{so}(2n,\mathbb{H}))\mid A, B\in M_n(\mathbb{R})\}$.  Therefore $\mathfrak{p}=\mathfrak{p}_0\otimes_\mathbb{R}\mathbb{C}\cong 
\{\left (\begin{smallmatrix}Z & W \\W & -Z\end{smallmatrix}\right)\in \mathfrak{so}(2n,\mathbb{C})\}$.  
Note that multiplication on the right by $J$ yields a complex structure on $\mathfrak{p}_0$.


\subsection{$F$-structures on $\SO(n,\mathbb{H})$} \label{f-structures}
Our aim in this section is to construct linear algebraic groups defined over number fields $F\subset \mathbb{R}$ 
such that the $\mathbb{R}$-points is isomorphic, as a Lie group, to $\SO(n,\mathbb{H})$.  This is most conveniently 
achieved by starting with quaternion algebras over number fields which do not split over $\mathbb{R}$.

Let $F$ be a subfield of $\mathbb{R}$ and let $\alpha, \beta\in F^\times$.
We denote by $\mathbb{D}:=\mathbb{H}^{\alpha, \beta}_F$ the quaternion algebra over $F$ generated by $i, j$ where $i^2=\alpha, j^2=\beta, ij=-ji=:k,$ so that $k^2=-\alpha\beta$.   
Sending $i$ to $i$ and $j\to k/\alpha$ yields an $F$-algebra isomorphism $\mathbb{H}_F^{\alpha, \beta}\to \mathbb{H}_F^{\alpha, -\beta/\alpha}$.  Also $\mathbb{H}^{\alpha,\beta}_F\cong\mathbb{H}^{\beta,\alpha}_F$.  In view of this, {\it we assume, without loss of generality, that $\beta<0$.}

One has an involutive anti-automorphism $\tau_c$ on $\mathbb{D}$ sending $q=q_0+iq_1+jq_2+kq_3$ to $\bar{q}=q_0-iq_1-jq_2-kq_3$.   
We have also the norm $N(q)=q\bar{q}=q_0^2-\alpha q_1^2-\beta q_2^2+\alpha\beta q_3^2$.    
The $F$-linear map $\tau_r:\mathbb{D}\to \mathbb{D}$ is the {\it reversion} $\tau_r(q)=q_0+iq_1-jq_2+kq_3$.

If $E\subset \mathbb{C}$ is an extension field of $F$ we denote $\mathbb{D}\otimes _FE=\mathbb{H}^{\alpha,\beta}_E$ by $\mathbb{D}_E$. 
Note that $\mathbb{D}$ splits over any field $E$ that contains  $F(a)\subset \mathbb{C}$ where $a^2=\alpha$.  In particular, $\mathbb{D}_\mathbb{R}\cong M_2(\mathbb{R})$ if $\alpha>0$. If both $\alpha$ and  $\beta$ are negative, then $\mathbb{D}_E$ is a division algebra 
over any extension field $E\subset \mathbb{R}$ of $F$ and we have an isomorphism $\mathbb{D}_E=\mathbb{H}^{\alpha,\beta}_E\cong \mathbb{H}^{a\alpha,b\beta}_E $ where $a, b\in (E^\times)^2$.

Suppose that $E\subset \mathbb{C}$ contains $F(a,b)$ where $a^2=\alpha, b^2=-\beta$, then 
$i\mapsto \left(\begin{smallmatrix}a& 0\\0 & -a\end{smallmatrix}\right), j\mapsto \left(\begin{smallmatrix}0& b\\-b & 0\end{smallmatrix}\right), $ defines an $E$-algebra 
isomorphism $\mathbb{D}_E\to M_2(E)$. 
Thus we have 
an $F$-algebra embedding $\psi_E:\mathbb{D}\hookrightarrow M_2(E)$.  
The embedding $\psi_E$ defines 
an $F$-algebra embedding (denoted by the same symbol) $\psi_E:M_n(\mathbb{D})\to M_{2n}(\mathbb{E})$ where $Q=Q_0+Q_1i+Q_2j+Q_3k\mapsto     \left ( \begin{matrix}Q_0+aQ_1 & bQ_2+abQ_3 \\-bQ_2+abQ_3 & Q_0-aQ_1\end{matrix}\right)$.     (As $\psi_E$ depends 
not only on $F$ but also on the choice of square roots of $\alpha, -\beta$, the notation is somewhat imprecise.) 
Again reversion on $\mathbb{D}$ induces an involutive anti-automorphism $Q\mapsto {}^t\tau_r(Q)$ of the $F$-algebra $M_n(\mathbb{D}).$  
One has the relation $\psi_E({}^t\tau_r(Q))={}^t\psi_E(Q)$ for all $Q$ in $M_n(\mathbb{D})$.    A matrix $A\in M_n(\mathbb{D})$ is said to be $\tau_r$-{\it hermitian} if ${}^t\tau_r(A)=A$.  
Writing $A=A_0+iA_1+jA_2+kA_3$, we see that $A$ is $\tau_r$-hermitian if and only if $A_2=-{}^tA,$ ${}^tA_r=A_r, r=0,1,3$.  If $A$ is $\tau_r$-hermitian 
and $\alpha>0$, then $\psi_E(A)$ is {\it real} symmetric.  (Recall that $\beta<0$ 
by hypothesis and so $b$ is real whereas $a$ is real if $\alpha>0$ and is purely imaginary if $\alpha<0$.)   

Let $A\in M_n(\mathbb{D})$ be $\tau_r$-hermitian and invertible.  The group $\SO(A,\mathbb{D}):=\{X\in M_n(\mathbb{D})\mid {}^t\tau_r(X)AX=A, ~\det(X)=1\}$ is the $F$-points of an $F$-algebraic group $\mathbf{G}$ such that $G(\mathbb{C})\cong \SO(2n,\mathbb{C})$. It is known that all $F$-structures on $SO(n,\mathbb{C})$ 
arise this way when $n>4$.  (See \cite[Ch. 18]{morris}.)
When $\alpha>0$ we have $G(\mathbb{R})=SO(A,\mathbb{D}_\mathbb{R})
\cong SO(p,2n-p)$ where the symmetric matrix $\psi_\mathbb{R}(A)\in M_{2n}(\mathbb{R})$ has exactly $p$ positive eigenvalues. In particular if $\psi_\mathbb{R}(A)$ is positive or negative definite, $G(\mathbb{R})$ is compact.  When $\alpha<0$ (and since $\beta<0$ by assumption) 
we have $\mathbb{D}_\mathbb{R}\cong \mathbb{H}$ and the 
group $G(\mathbb{R})$ is isomorphic to the type AIII real form 
$SO(n,\mathbb{H})\cong \SO^*(2n)$ of $SO(2n,\mathbb{C})$.

Suppose that $P\in M_n(\mathbb{D})$ is a non-singular matrix such that ${}^t\tau_r(P)AP=\lambda A$ for some 
$\lambda\in F$.  
Then it is easily verified that $Q\mapsto PQP^{-1}$ defines an $F$-automorphism of $\SO(A, \mathbb{D})$. 

 
 \subsection{$F$-rational Cartan involution} \label{frationalcartan}
Recall that $J=jI_n\in M_n(\mathbb{H})$ and that 
$Q\mapsto JQJ^{-1}$ is the Cartan involution $\theta$ of $\SO(n,\mathbb{H})$ that fixes the 
standard maximal compact subgroup $K=\{Z+Wj\in \SO(n,\mathbb{H})\mid Z,W\in M_n(\mathbb{R})\}.$ 
 Assume that $\alpha,\beta\in F$ are negative and let $\mathbb{D}=\mathbb{H}^{\alpha,\beta}_F$.  It is readily seen that $\theta$ 
is also the Cartan involution of $\SO(I_n,\mathbb{D}_\mathbb{R})$ that fixes the standard maximal 
compact subgroup $K=\{Z+Wj\in \SO(n,\mathbb{D}_\mathbb{R})\mid Z, W\in M_n(\mathbb{R})\}.$  In fact $\theta$ is induced from 
an $F$-rational involution $\theta:\SO(I_n,\mathbb{D})\to \SO(I_n,\mathbb{D})$ which fixes the subgroup 
$K_F:=\{ Z+Wj\in \SO(I_n,\mathbb{D})\mid 
Z,W\in M_n(F)\}.$ The group $K_F$ is an $F$-algebraic group whose $\mathbb{R}$-points  equals $K$.

Let $A={}^t\tau_r(P).P$ where $P\in M_n(\mathbb{D})$ is non-singular.  Then $A$ is $\tau_r$-hermitian and non-singular.   
One has an $F$-isomorphism $\SO(A,\mathbb{D})\to \SO(I_n,\mathbb{D})$ defined as $Q\mapsto P Q P^{-1}$.   
It follows that $K_A:=P^{-1} K_FP$ is an $F$-algebraic group whose $\mathbb{R}$-points is the maximal compact 
subgroup $P^{-1}KP\subset \SO(A, \mathbb{D}_\mathbb{R})$.  
We shall prove the existence 
of an $F$-rational Cartan involution for a more general class of $\tau_r$-hermitian matrices.  This will be 
preceded by some preliminary observations.

\begin{lemma}\label{diagonal-A}
Let $\alpha, \beta\in F^\times$ with $\beta<0$, $\mathbb{D}=\mathbb{H}_F^{\alpha,\beta}$ and let $A$ be any $\tau_r$-hermitian matrix.  There exists a  non-singular matrix $P\in M_n(\mathbb{D})$ and a 
diagonal matrix $D=\diag(d_1,\ldots, d_n) \in M_n(\mathbb{D})$ such that ${}^t\tau_r(P).A.P=D$.  
\end{lemma}

\begin{proof} 
Regard $V=\mathbb{D}^n$ as a right $\mathbb{D}$-module.   One has a $\tau_r$-{\it hermitian} pairing $\langle.,.\rangle:V\times V\to \mathbb{D}$ defined as $\langle u,v\rangle :={}^t\tau_r(u) A v$.  It is $\mathbb{D}$-linear in the second argument and $\tau_r$-twisted $\mathbb{D}$-linear in the first, that is, 
$\langle uq,vq'\rangle=\tau_r(q) \langle u,v\rangle q'$ for $q,q'\in \mathbb{D}$.  If $e_1,\ldots, e_n\in F^n\subset \mathbb{D}^n$ 
denotes the standard basis of $V$ over $\mathbb{D}$, then $A$ equals the matrix $(\langle e_i,e_j\rangle)$.
Since $A$ is non-singular, the pairing is non-degenerate.  
It is easy to show, using a Gram-Schmid orthogonalisation argument starting with the standard basis,  the existence of an orthogonal basis $\mathcal{B}:=\{v_1,\ldots, v_n\}$ for $V$ (as a right $\mathbb{D}$-module), although it is not 
in general possible to obtain an ortho{\it normal} basis. 
Let  $d_l:=\langle v_l,v_l\rangle={}^t\tau_r(v_l)Av_l$.    
We set $P\in M_n(\mathbb{D})$ to have columns the vectors $v_1,\ldots, v_n$. 
Then ${}^t\tau_r(P)AP=\diag(d_1,\ldots,d_n)=:D$. 
\end{proof}

The $v_l$---and hence the matrix $P$---is uniquely determined by the orthogonalisation process.  Inductively, 
the $v_l$ are determined by the following requirements:
$v_1=e_1$ and $v_l\perp v_k, k<l, v_l-e_l$ belongs to the $\mathbb{D}$-span of $e_1,\ldots, e_{l-1}$.

The group $\SO(A, \mathbb{D})$ is $F$-isomorphic to $\SO(D, \mathbb{D})$ 
where $D$ is as in the above lemma.  An $F$-isomorphism $\SO(A,\mathbb{D})\to \SO(D,\mathbb{D})$  is given by $Q\mapsto P^{-1}QP$. 

\begin{lemma}  \label{diagonal-P}
Suppose that $\alpha, \beta<0$.
Let $D=\diag(d_1,\ldots, d_n)\in M_n(\mathbb{D})$ be non-singular and $\tau_r$-hermitian.  
Then there exists a $\tau_r$-hermitian {\it diagonal} matrix $P\in M_n(\mathbb{D}_\mathbb{R})$ such that $D={}^t\tau_r(P)P=P^2$. 
\end{lemma}
\begin{proof}  It suffice to show that any $d\in \mathbb{D}$ such that $\tau_r(d)=d=a+bi+ck$ is expressible as $d=\tau_r(u).u=u^2$ 
with  $u=\tau_r(u)\in \mathbb{D}_\mathbb{R}$. 
We need to find $u=x+yi+wk$ such that $d=(x+yi+wk)(x+yi+wk)$.  We may (and do) assume that $bc\ne 0$.  This is equivalent to solving the following 
system of equations over $\mathbb{R}$:
\[ x^2+\alpha y^2-\alpha\beta w^2=a,\]
\[2xy=b,\]
\[2xw=c.\]
We obtain $2w=c/x, 2y=b/x$ and so substituting in the first equation implies that 
$f(x):=4x^2+\alpha b^2/x^2-\alpha\beta c^2/x^2=a$.  Since $\alpha, -\alpha\beta<0$, the function $f:\mathbb{R}_{>0}\to 
\mathbb{R}$ is continuous and onto.  Let $f(x_0)=a$.  Then $u=x_0+bi/(2x_0)+ck/(2x_0)$ 
is a solution.   
\end{proof}

\begin{remark} {\em
The value of $x_0\in \mathbb{R}_{>0}$ in the above proof is unique.  This means that $u$ is unique up to sign.  The 
latter statement is evidently valid if $bc=0$.  It follows that 
each diagonal entry of $P$ uniquely determined up to a sign $\pm$. }
\end{remark}

\begin{theorem} \label{f-cartan}
Let $A$ be any $\tau_r$-hermitian non-singular matrix in $M_n(\mathbb{D})$ where $\mathbb{D}=\mathbb{H}_F^{\alpha,\beta}$ with 
$\alpha, \beta\in F$ both negative.  Suppose that ${}^t\tau_r(P_0)AP_0=D=\diag(d_1,\ldots, d_n)$ where 
$N(d_l)\equiv N(d_1)\mod (F^\times)^2~\forall l\le n$.  Then there exists an $F$-involution $\theta_A$ of $\SO(A,\mathbb{D})$ which defines 
a Cartan involution of $\SO(A,\mathbb{D}_\mathbb{R})$.  In fact $\theta_A$ is the restriction to $\SO(A,\mathbb{D})$ of 
$\iota_Y$ where $Y=cPJP^{-1}\in M_n(\mathbb{D})$ for a suitable  $P\in M_n(\mathbb{D}_\mathbb{R})$ and $c\in \mathbb{R}$.
\end{theorem}
\begin{proof}  Since $\SO(A, \mathbb{D})$ is $F$-isomorphic to $\SO(D,\mathbb{D})$ it suffices to show that 
$\SO(D,\mathbb{D})$ admits an $F$-rational Cartan involution $\theta_D$ with $\Fix(\theta_D)=K_D$. 
By Lemma \ref{diagonal-P}, there exists diagonal matrix $P=\diag(p_1, \ldots, p_r)\in M_n(\mathbb{D}_\mathbb{R})$ 
such that ${}^t\tau_r(P)
=P, D=P^2={}^t\tau_r(P)P.$   In particular, since $\iota_J$ is the Cartan involution that fixes $K\subset \SO(I_n,\mathbb{H}_\mathbb{R}^{\alpha,\beta}),$  it follows that  $\iota_{P^{-1}JP}$ is the Cartan involution that fixes $P^{-1}KP$. 
It remains to find a $c\in \mathbb{R}$ so that $\iota_Y=\iota_{P^{-1}JP}$ is an 
$F$-involution where $Y:=cP^{-1}JP\in M_n(\mathbb{D})$. 
By our hypothesis on $D$, there exist positive numbers $t_l\in F^\times, 1\le l\le n,$ with $t_1=1, t_l>0~
\forall l$ such 
that $t_l^2N(d_l)=N(d_1)~\forall l.$ 
Let $c_l=t_lN(p_l)\in F^\times, 1\le l\le n.$  Then $c_l^2=t_l^{2}N(p_l)^2= t_l^{2}N(d_l)=N(d_1)=c_1^2$.  As $c_l>0$, 
we have $c_l=c_1=:c$ for all $l\le n$.  It remains to show that the diagonal matrix $Y=cP^{-1}JP$ is 
in $M_n(\mathbb{D})$.  We compute $l$-th diagonal entry---call it $y_l$---of $Y$ using $p_l=\tau_r(p_l)=j\tau_c(p_l)j/\beta$.  We have $y_l=cp_l^{-1} j p_l=c \frac{\tau_c(p_l)}{N(p_l)}jp_l
=\frac{cj} {N(p_l)\beta}.j\tau_c(p_l)jp_l 
=jt_lp_l^2=jt_ld_l\in \mathbb{D}$.   Thus we may take $\theta_A$ to be the 
restriction of $\iota_Y$.
\end{proof}

\subsection{Commuting pairs of $F$-Involutions} \label{commutinginvolutions}
We keep the notations of \S2.2.
We shall now construct $F$-involutions of $G=\SO(A,\mathbb{D})$ where $\mathbb{D}=\mathbb{H}_F^{\alpha,\beta}$, 
with $\alpha,\beta<0$, and $A\in M_n(\mathbb{D})$ is $\tau_r$-hermitian.   
The $F$-involutions which commute with the $F$-rational Cartan involution $\theta_A$ would play 
a crucial role in the sequel.  

Let $A\in M_n(\mathbb{D})$ be $\tau_r$-hermitian and non-singular.  In view of Lemma \ref{diagonal-A}, we may 
assume, without loss of generality, that 
$A=\diag(a_1,\ldots, a_n)\in M_n(\mathbb{D})$ is diagonal.  Then $a_l=x_l+y_li+z_lk~\forall l.$ It is readily checked that for $S:=\diag(\epsilon_1,\ldots, \epsilon_n), \epsilon_l\in \{-1,1\}, $ 
satisfies the conditions $S^2=I, {}^t\tau_r(S)=S,~A={}^t\tau_r(S)AS=SAS$ and so $Q\mapsto 
S^{-1}QS=SQS$ defines an involutive $F$-automorphism of $\SO(A, \mathbb{D})$.   These involutions 
will be referred to as the {\it sign} involutions. 

More generally, $Q\mapsto S^{-1}QS$ 
is an involutive $F$-automorphism if $S\in M_n(\mathbb{D})$ satisfies the following conditions: \\
\[S^2=-\lambda I_n, ~\textrm{for some~}\lambda\in F^\times  \hfill \eqno(1)\]
\[{}^t\tau_r (S)AS=\mu A, \textrm{~for some~}\mu \in F^\times.\hfill \eqno(2)\]

First we shall classify all such $S=\diag(s_1,\ldots, s_n)$ that are diagonal.  By (1), $s_k^2=-\lambda, \forall k$.  
If $s_k \in F^\times$ for some $k$, then $\lambda<0$. So  for any $l$, $s_l^2>0$ and so $s_l\in F^\times $.  Thus $s_i=\pm s_j ~\forall i,j$ and so that $S=\diag (\epsilon_1s,\ldots, \epsilon_ns)$ where $s=|s_1|$ and  
each $\epsilon_l$ equals $1$ or $-1$.   The involution defined by $S$ is therefore a sign involution.   

Suppose that $s_k\notin F^\times$ for some $k$.  Then $s_k^2=-\lambda\in F^\times$ implies that $s_k\in Fi+Fj+Fk$ and 
$\lambda=-s_k^2=N(s_k)>0$.  Therefore for {\it any} $l\le n$, $s_l^2<0$ and so $s_l\in Fi+Fj+Fk,  N(s_l)=\lambda$.  Now equation (2) says that 
$\tau_r(s_k) a_k s_k=\mu a_k$ and so $\mu^2.N(a_k)=N(\tau_r(s_k) a_k s_k)=N(a_k)N(s_k)^2=\lambda^2N(a_k) $ and 
so $\mu=\pm \lambda$.   Using $\tau_r(s)=j^{-1}\tau_c(s)j=-j^{-1}sj ~\forall s\in Fi+Fj+Fk,$ we obtain that $\mu a_k=
\tau_r(s_k)a_ks_k=-j^{-1}s_kja_ks_k$ which implies $\mu s_k. ja_k=  s_k. (-s_k).ja_k s_k=\lambda .ja_k. s_k=\pm \mu ja_k. s_k$, that is, either $s_k.ja_k=ja_k.s_k ~\forall k$ in which case $\lambda=\mu$, or,  $s_k.ja_k=-ja_k.s_k~\forall k$ in which case $\lambda=-\mu$.  Equivalently $S.jA=jA.S $ or $S.jA=-jA.S$.

 Since $s_k\in Fi+Fj+Fk, $ if $s_k$ commutes with $ja_k$ then $s_k\in F(ja_k)$. In this case $S=DjA$ 
where $D\in M_n(F)$ is a diagonal matrix. 

 It is clear that $s_k$ anti-commutes with $ja_k$ if and only if 
$s_k\perp (ja_k) $ in $Fi+Fj+Fk\subset \mathbb{D}$ 
with respect to the inner product on the $F$-vector space $\mathbb{D}$, which  
is  associated to the quadratic form 
$q\mapsto N(q)$, which is positive definite since $\alpha,\beta<0$. 
This is immediate from comparing the coefficient of $1$ on both sides of $s_k (ja_k)=-(ja_k)s_k$; observe that ${}^t\tau_r(A)=A$ implies that 
$ja_l\in Fi+Fj+Fk$ for all $l$. 
The involution on $\SO(A,\mathbb{D})$ induced by $S$ will be called an {\it involution 
of even type} if 
$SjA=jAS$, equivalently $\lambda=\mu$.  If $SjA=-jAS$, then the involution induced by $S$ will be referred to as an 
{\it involution of odd type}, equivalently $\lambda=-\mu$.

\begin{lemma} \label{involutiontype}
Let $A=\diag(a_1,\ldots, a_n)\in M_n(\mathbb{H}_F^{\alpha, \beta}), \alpha, \beta<0$ be a non-singular diagonal $\tau_r$-hermitian matrix.  
(i) An $F$-involution of $\SO(A,\mathbb{D})$ of even type exists if and only if $N(a_j)\equiv N(a_k) \mod 
(F^\times)^2$ for all $1\le j,k\le n$.  \\
 (ii)  An involution of odd type exists if $N(a_j)\equiv N(a_k)\mod (F^\times)^2$ for all $j, k$. \\
(iii) Any two $F$-involutions of even type (possibly for different values of $\lambda$) commute. 
Also any $F$-involution of even type commutes with any $F$-involution of odd type.
\end{lemma}
\begin{proof}
(i) Let $N(a_k)=t_k^2N(a_1), k\le n$ for suitable $t_k\in F^\times$.  
Set $S=\diag(s_1,\ldots, s_n)$ with $s_k:=t_k^{-1}ja_k~\forall k\le n$.    Then $S$ satisfies equations (1) and (2) 
with $\lambda = N(ja_1)=\mu $ and so $S$ induces an involution of even type.  
The discussion preceding the statement of the lemma establishes the converse part.

(ii) Suppose that $N(a_k)=t_k^2N(a_1), 1\le k\le n.$  Then $N(t_k^{-1}ja_k)=N(ja_1)$ for all $k\le n$.  
Thus $ja_k\mapsto  t_kja_1$ defines an isomorphism of inner product spaces $Fja_k\to Fja_1.$   It follows 
their orthogonal complements  $(Fja_1)^\perp, (Fja_k)^\perp$ in $Fi+Fj+Fk$ are isometric, that is, they are isomorphic as inner product spaces. (See \cite[Ch. XV, Theorem 10.2]{lang}.)  Fix linear isometries 
$f_k:(Fja_1)^\perp \to (Fja_k)^\perp$  and choose 
$s_1\in (Fja_1)^\perp $ to be any non-zero element.  If $s_k:=f_k(s_1), 2\le k\le n$, then $S:=\diag(s_1,\ldots, s_n)$ 
satisfies (1) and (2) with $-\mu=\lambda =N(s_1)$, and so $S$ induces an involution of odd type.  

(iii) If $S$ determines an $F$-involution $\sigma$ of even type and $S'$ an $F$-involution $\sigma'$ of odd type, 
then, by the discussion preceding the statement of the lemma, we have a diagonal matrix $D\in M_n(F)$ such that 
$S=DjA$.  Since $D$ commutes with any diagonal matrix in $M_n(\mathbb{D})$ we have 
$SS'=DjAS'=-DS'jA=-S'DjA=-S'S$.  It follows that $\iota_S\iota_{S'}=\iota_{SS'}=\iota_{-S'S}=\iota_{S'S}=\iota_S\iota_{S'}$ and so $\sigma$ and $\sigma'$ commute.  
We leave to the reader the verification that any two involutions of even type commute.  
\end{proof}

\begin{remark}\label{cartaniseven}
{\em
(i) 
We observe that, assuming $A=\diag(a_1,\ldots,a_n)$, the Cartan involution $\theta_A$ constructed in Theorem \ref{f-cartan} is of even type.  This is because, in the notation of that theorem, $\theta_A$ is the restriction of $\iota_Y$ to $\SO(A,\mathbb{D})$ where $Y=\diag(y_1,\ldots, y_l)$ with $y_l=t_ljd_l\in Fjd_l~\forall l\le n$.  The value of 
$\lambda$ for which equation (1) satisfied by $\theta_A$ is obtained as $\lambda=N(y_1)=-\beta N(d_1)$.  In particular 
any involution of even or odd type commutes with $\theta_A$.

(ii)   Suppose that $\sigma$ is an involution of even type induced by a diagonal matrix $S$ that 
satisfies equations (1) and (2).   Then 
$S=CY$ where $C\in M_n(F)$ is a diagonal matrix of the form $\diag (\epsilon_1c,\ldots,\epsilon_nc),$ 
$\epsilon_k\in \{1,-1\}$.  It follows that $\sigma=\varepsilon\circ \theta_A$ where $\varepsilon$ 
is a sign involution.
}
\end{remark}

\subsection{Fixed points of involutions}\label{fixedpoints}
Let $F$ be a real number field.
Let $A\in M_n(\mathbb{D})$ be $\tau_r$-hermitian and non-singular, where $\mathbb{D}=\mathbb{H}_F^{\alpha,\beta}$, $\alpha, \beta\in F$ both negative. 
Without loss of generality we may (and do) assume that $A$ is diagonal, say, $A=\diag(a_1,\ldots, a_n)$; see \S \ref{frationalcartan}.   Assume that $N(a_l)\equiv N(a_1) \mod (F^\times)^2$.  Let $S\in M_n(\mathbb{D})$ be non-singular 
diagonal matrix $\diag(s_1,\ldots, s_n)$ that satisfies (1), (2) of \S \ref{commutinginvolutions}.   Denote by $\sigma$ the 
involution on $\mathbf{G}:=\SO(A,\mathbb{D})$ defined by conjugation by $S$. 

The group $C:=\langle \sigma\rangle$ is the $F$-group 
$\textrm{Spec}(F[t]/\langle t^2-1\rangle)$.   
 We denote 
by $\mathbf{G}(\sigma):=\textrm{Fix}(\sigma)$, the fixed subgroup of $\mathbf{G}$ fixed by the action of $C$ on $\mathbf{G}$. 
It is an $F$-algebraic subgroup of $\mathbf{G}$ whose $F$-algebra of regular functions is the algebra of co\"{i}nvariants $F[G]/\langle f-\sigma(f)\mid f\in F[G]\rangle$ for the action of $C$ on $F[G]$.

The $\mathbb{R}$-points $G(\sigma)(\mathbb{R})$ of 
the $F$-group $\mathbf{G}(\sigma)$ is the fixed subgroup  
of $G(\mathbb{R})$ for the action of $C_\mathbb{R}=
\text{Spec}(\mathbb{R}[t]/\langle t^2-1\rangle)$ on $G(\mathbb{R})$.  That is,  
$G(\sigma)(\mathbb{R})=\{Q\in G(\mathbb{R})\mid 
\sigma(Q)=Q\}$. 
 
For our applications it is important to know if the  group $G(\sigma)(\mathbb{R})$ acts by preserving 
the orientation on the symmetric space $X(\sigma):=G(\sigma)(\mathbb{R})/K(\sigma)$, where $K(\sigma)=K\cap G(\sigma)(\mathbb{R})$.  Of course it is enough to check this for a set of elements belonging 
to each connected component of $G(\sigma)(\mathbb{R})$.  The requirement that $G(\sigma)(\mathbb{R})$ act 
by preserving the orientation on $X(\sigma)$ is called condition {\it Or} in \cite[Theorem 4.11]{rs}.

We first consider $S$ so that $\sigma $ is a sign involution.  If suffices to consider the case when $S=S_l:=\left(\begin{smallmatrix} I_l&\\&-I_{n-l}\end{smallmatrix}\right)$.  In this case the $F$-points of $G(\sigma)$ is easily 
determined: An element $Q\in \SO(A, \mathbb{D})$ is fixed by $\sigma$ if and only if $Q$ is block diagonal 
$\left(\begin{smallmatrix} X&\\&Y\end{smallmatrix}\right)$ where $X\in \SO(A_l,\mathbb{D}), Y\in \SO(A'_{n-l}, \mathbb{D})$. 
Here $A_l=\diag(a_1,\ldots, a_l), A'_{n-l}=\diag(a_{l+1},\ldots, a_n)$.   
In this case the group $G(\sigma)(\mathbb{R})$ is the connected group $\SO(A_l,\mathbb{D}_\mathbb{R})\times \SO(A'_{n-l},\mathbb{D}_\mathbb{R})\cong \SO(l,\mathbb{H})\times \SO(n-l,\mathbb{H}), $ and 
$X(\sigma)\cong \SO(l,\mathbb{H})/U(l)\times \SO(n-l,\mathbb{H})/U(n-l)$.  
Since $G(\sigma)(\mathbb{R})$ is connected, the condition {\it Or} is satisfied.  We remark that $X(\sigma)$ is 
hermitian symmetric and in fact the inclusion $X(\sigma)\hookrightarrow X=G(\mathbb{R})/K=\SO(A,\mathbb{D}_\mathbb{R})/U(n)$ is complex analytic.  This is because $\sigma(z)=z$ for all $z\in Z_K\cong \mathbb{S}^1$ the centre of 
$K=U(n)$.  (Cf. \cite[Remark 4.8(ii)]{rs}.) 

Next we turn to $F$-involutions of the even type.   In this case, in view of Remark \ref{cartaniseven}, 
we see that any such involution  $\sigma$ equals $\varepsilon\theta_A$ where $\varepsilon$ is a sign 
involution.  For notational convenience we assume that $\varepsilon$ is induced by $S_l=\left(\begin{smallmatrix}
 I_l&\\&-I_{n-l}\end{smallmatrix}\right)$.   
 
The following observation reduces computation 
of $G(\sigma)(\mathbb{R})$ to the special case 
when $A=I$ and $\alpha=\beta=-1$.    

In view of Lemma \ref{diagonal-P}, we have an isomorphism of Lie groups:
$\SO(A,\mathbb{D}_\mathbb{R})\to \SO(I_n, \mathbb{D}_\mathbb{R})$
given by $Q\mapsto PQP^{-1}$ where $P$ is an 
invertible $\tau_r$-hermitian in $\SO(A, \mathbb{D})$ 
such that (the diagonal matrix) $A={}^t\tau_r(P).P=P^2$.  Also this isomorphism is $\mathbb{Z}/2\mathbb{Z}$-equivariant where the action on $\SO(A,\mathbb{D}_\mathbb{R})$ is given by $\sigma$ and on $\SO(I_n,\mathbb{D}_\mathbb{R})$ 
by $\varepsilon \theta_{I_n}$.  
This is because $\theta_A$ is the conjugation by $cP^{-1}JP$ with $c\in \mathbb{F}^\times$
(see Theorem \ref{f-cartan}), and $S_l$ commutes with $P$ (since $P$ is diagonal).  Therefore $G(\sigma)(\mathbb{R})\to \text{Fix}(\varepsilon\theta_{I_n}), Q\mapsto PQP^{-1}$ is an isomorphism.  Finally $\SO(I_n,\mathbb{D}_\mathbb{R})\cong 
\SO(I_n,\mathbb{H}_\mathbb{R}^{-1,-1})=\SO(n,\mathbb{H})$ where $i\mapsto ai, j\mapsto bj, a,b\in \mathbb{R}_{>0}$ with $a^2=-\alpha, 
b^2=-\beta$.  Clearly this isomorphism is again equivariant with respect to the action of
$\langle \varepsilon\theta_{I_n}\rangle$ on the domain and that of $\langle\varepsilon\theta\rangle$ the target. (Recall from \S \ref{realform} that $\theta$ is the Cartan involution that fixes the standard maximal compact 
subgroup $K=\{X+Zj\in \SO(n,\mathbb{H})\mid X,Z \in M_n(\mathbb{R})\}$.) 
Hence $G(\sigma)(\mathbb{R})\cong \textrm{Fix}(\varepsilon\theta)$.

We now compute $\textrm{Fix}(\varepsilon\theta)$.
Write $X\in M_n(\mathbb{R})$ as a block matrix $\left(\begin{smallmatrix} X_{11} &X_{12}\\ 
X_{21} &X_{22}\end{smallmatrix}\right)$ where $X_{11}\in M_l(\mathbb{R}), X_{22}\in M_{n-l}(\mathbb{R})$.  
Then $\varepsilon (X)=S_lXS_l^{-1}=\left(\begin{smallmatrix} X_{11} &-X_{12}\\ 
-X_{21} &X_{22}\end{smallmatrix}\right)$.  Thus $X=\varepsilon( X)$ if and only if $X_{12}=0, X_{21}=0$ and $X=-\varepsilon (X)$ if and only if $X_{11}=0, X_{22}=0$.  

Now let $Q=X+Yi+Zj+Wk\in \SO(n,\mathbb{H})$ with $X,Y,Z,W\in M_n(\mathbb{R})$.  Then 
$\theta(Q)=X-Yi+Zj-Wk$.  So $\varepsilon(\theta(Q))=Q$ if and only if, writing $X, Y, Z, W\in M_n(\mathbb{R})$ as block 
matrices as in the previous paragraph, we have $X_{12}=Z_{12}=0, X_{21}=Z_{21}=0, Y_{11}=W_{11}=0, 
Y_{22}=W_{22}=0$. Hence $K(\varepsilon\theta):=\textrm{Fix}(\varepsilon\theta)\cap K=\{X+Zj\in \SO(n,\mathbb{H})\mid X_{12}=Z_{12}=0, X_{21}=0=Z_{21}\}
\cong U(l)\times U(n-l)$ since $X$ is skew symmetric and $Z$ is symmetric as $X+Zj\in \SO(n,\mathbb{H})$.    
Since  $K(\varepsilon\theta)$ is connected and is a maximal compact subgroup of $\textrm{Fix}(\varepsilon\theta)$, the latter group is connected.

We claim that $\textrm{Fix}(\varepsilon\theta)=U(l,n-l)$.  We need only show that the Lie algebra of $\text{Fix}(\varepsilon\theta)$ is isomorphic to $\mathfrak{u}(l,n-l) $.  
To establish the claim, consider the Cartan decomposition $\mathfrak{so}(n,\mathbb{H})
=\mathfrak{k}_0\oplus \mathfrak{p}_0,$ defined by $\theta$.  
As $\varepsilon$ commutes with $\theta$, the spaces 
$\mathfrak{k}_0,\mathfrak{p}_0\subset \mathfrak{so}(n,\mathbb{H})$ are stable by $\varepsilon$, and,   
$\text{Lie}(\textrm{Fix}(\varepsilon\theta))\subset \mathfrak{so}(n,\mathbb{H})$ decomposes as a direct sum of 
$(+1)$-eigenspaces of $\varepsilon|_{\mathfrak{k}_0}$ and the $(-1)$-eignspace of $\varepsilon|_{\mathfrak{p}_0}$.   This 
is in fact the Cartan decomposition of $\text{Lie}(\textrm{Fix}(\varepsilon\theta))$. 
Denote these subspaces by $\mathfrak{k}_0(\epsilon\theta), \mathfrak{p}_0(\epsilon\theta)$ respectively.  
A direct calculation shows that $\mathfrak{k}_0(\epsilon\theta)$ consists of 
$X+Zj\in \mathfrak{so}(n,\mathbb{H})$ with $X=\left(\begin{smallmatrix} X_{11}&\\ 
& X_{22} \end{smallmatrix}\right)$ skew symmetric  and  $Z=\left(\begin{smallmatrix} Z_{11}&\\ 
& Z_{22} \end{smallmatrix}\right)$ symmetric, (where diagonal blocks have sizes $l$ and $n-l$).    
Similarly $\mathfrak{p}_0(\epsilon\theta)$ consists of $Yi+Wk\in \frak{so}(I_n, \mathbb{D}_\mathbb{R})$ 
with $Y=\left(\begin{smallmatrix} &Y_{12}\\  Y_{21}&  \end{smallmatrix}\right), W=\left(\begin{smallmatrix} &W_{12}\\ 
W_{21}& \end{smallmatrix}\right)$; thus $Y,W$ are real skew symmetric matrices.    
Now it is readily verified that  the Lie algebra of $\textrm{Fix}(\epsilon\theta)$ is isomorphic to $\mathfrak{u}(l,n-l)$.
Thus $\mathfrak{g}_0(\sigma)\cong \mathfrak{u}(l,n-l)$.   

It remains to consider $\textrm{Fix}(\sigma)$ when $\sigma$ is an involution of odd type.  
As in the case of involutions of even type, it suffices to consider the case when $A=I_n, \alpha=\beta=-1$, 
$\sigma$ is induced by conjugation by a diagonal matrix $S=\diag(s_1,\ldots, s_n) \in M_n(\mathbb{H})$ such that $S^2=-I_n$, 
$SJ=-JS$ where $J:=jI_n$ and $N(s_l)=1~\forall l$.  Write $s_l=ie(\lambda_l)\in \mathbb{R}i+\mathbb{R}k$, 
where $e(\lambda):=\cos(\lambda)+j\sin(\lambda)\in \mathbb{R}[j]$.   

Suppose that $Q=(q_{lm})\in M_n(\mathbb{H})$ is such that $\sigma(Q)=SQS^{-1}=Q$.  Then, writing $q_{lm}= z_{lm}+iw_{lm}$ with 
$z_{lm},w_{lm}\in \mathbb{R}[j]$, a straightforward computation yields, for $1\le l,m\le n$,  $z_{lm}=x_{lm}e(\frac{\lambda_m-\lambda_l}{2}), 
w_{lm}=y_{lm}e(\frac{\lambda_l+\lambda_m}{2})$ where $X=(x_{lm}), Y=(y_{lm}) \in M_n(\mathbb{R})$.  
Set $T:=\diag(e(\lambda_1/2), 
\cdots,e(\lambda_n/2))$.  Then we have 
$Q=T^{-1}XT+i TYT=T^{-1}XT+T^{-1}iYT=T^{-1}(X+iY)T$.  
Conversely, if $Q=T^{-1}(X+iY)T$ with $X,Y\in M_n(\mathbb{R})$, then $\sigma(Q)=Q$ using $S=iT^2=T^{-2}i$.  
Thus $V:=\{Q\in M_n(\mathbb{H})\mid \sigma(Q)=Q\}=T^{-1}M_n(\mathbb{R}[i])T$.   
Since ${}^t\tau_r(T)=T^{-1}$, it follows  that ${}^t\tau_r(Q)=T^{-1}.{}^t(X+iY).T$ for all $Q\in V$.  

We are ready to compute the group $G(\sigma)(\mathbb{R})\subset \SO(n,\mathbb{H})$.   
Clearly $G(\sigma)(\mathbb{R})=\SO(n,\mathbb{H})\cap V=\{T^{-1}CT\in V\mid {}^t\tau_r(T^{-1}CT).T^{-1}CT=I_n, C\in M_n(\mathbb{R}[i])\}
=\{T^{-1}CT\in V\mid {}^tC.C=I_n\}=T^{-1}O(n,\mathbb{R}[i])T\cong O(n,\mathbb{C})$.   We conclude that 
the group $G(\sigma)(\mathbb{R})$ has exactly two components.  The Cartan involution $\theta$ restricts 
to a Cartan involution on $G(\sigma)(\mathbb{R})$ and we see that $\mathfrak{p}_0(\sigma)
=\{T^{-1}iYT\mid {}^tY=-Y\in M_n(\mathbb{R})\}$.

\begin{lemma} \label{orforoddtype}
The action of the group $G(\sigma)(\mathbb{R})\cong O(n,\mathbb{C})$ on $X(\sigma)\cong O(n,\mathbb{C})/O(n)$ 
preserves the orientation if and only if $n$ is odd.  Thus involutions of odd type satisfies condition {\em Or}
of \cite[Theorem 4.11]{rs} precisely when $n$ is odd. 
\end{lemma}
\begin{proof}
It suffices to show that 
the action of the element $g=T^{-1}\diag(-1,1\ldots, 1)T\in G(\sigma)(\mathbb{R})$, given by the isotropy representation on 
$\mathfrak{p}_0(\sigma)=\{T^{-1}iYT\mid {}^tY=-Y\in M_n(\mathbb{R})\},$ 
is orientation reversing if and only if $n$ is even.  This is equivalent to showing that the element $\diag(-1,1,\ldots, 1)(=g)$ 
reverses the orientation on the space of all $n\times n$ real skew symmetric matrices if and only if $n$ is even. 
The lemma now follows by a straightforward computation. 
\end{proof}

\begin{remark}{\em 
When $\sigma$ is a sign involution or an involution of even type, $X(\sigma)$ is hermitian symmetric. 
On the other hand, if $\sigma$ is of odd type, then $X(\sigma)$ is a Lagrangian submanifold of $X=\SO(I_n,\mathbb{D}_\mathbb{R})/K$.
This is because the $G(\mathbb{R})$-invariant integrable almost complex structure $\mathfrak{J}$ on $X$ is obtained by conjugation by $Z:=(I_n-J)/\sqrt{2}=e(-\pi/4)I_n\in K$ on the tangent space at the identity coset, namely $\mathfrak{p}_0$.   If $Q=iQ_1+kQ_3\in \mathfrak{p}_0$, then $\mathfrak{J}(iQ_1+kQ_3)=e(-\pi/4)Qe(\pi/4)=Qe(\pi/2)=Q.J=-iQ_3+kQ_1$. So if $Q=T^{-1}iYT\in \mathfrak{p}_0(\sigma)$, as $Z$ commutes with $T$, we have $\mathfrak{J}(Q)
=T^{-1}iYT.J=T^{-1}kYT.$     
Since $iX$ and $kY$ are orthogonal for the 
Riemannian metric on $\mathfrak{p}_0$, and since conjugation by $T\in K$ is an 
isomtery of $\mathfrak{p}_0$, we conclude that $X(\sigma)$ is 
a Lagrangian submanifold of $X$. 
In view of the $G$-invariance of the Riemannian 
 metric, it follows that $C(\sigma)\subset X_\Lambda$ is a Lagrangian submanifold. 
}
\end{remark}


\subsection{Restriction of scalars and arithmetic lattices} \label{restrictionofscalars} 

Let $F\subset \mathbb{R}$ be a number field, i.e. a finite extension of $\mathbb{Q}$.   As usual $\bar{\mathbb{Q}}$ denotes the field of 
algebraic numbers.

We shall denote by $V_\infty$ the set of all real or complex embeddings of $F$ and by $S_\infty$ the set of all real embeddings.  The inclusion $F\hookrightarrow 
\mathbb{R}$ will be denoted by $\iota\in V_\infty$.   
We denote by $\mathcal{O}_F$ the ring of integers in $F$.  
Let $\mathbf{G}$  be a semisimple algebraic group defined over 
$F$, we denote by $G(\mathcal{O}_F)$ the arithmetic subgroup of $G(\mathbb{R})$ namely the $\mathcal{O}_F$-points of $\mathbf{G}$.    
For $\sigma\in V_\infty$, we denote by $\mathbf{G}^\sigma$ the corresponding algebraic group defined over $\sigma(F)$.   
Let $\mathcal{G}$ be the $\mathbb{Q}$-algebraic group got by restriction of scalars 
from $F$ to $\mathbb{Q}$, that is, $\mathcal{G}:=\res_{F/\mathbb{Q}}(\mathbf{G})$.  Then by a 
theorem of Borel and Harish-Chandra, 
$\mathcal{G}(\mathbb{Z})$, the group of $\mathbb{Z}$-points in $\mathcal{G}$, is a lattice 
in $\mathcal{G}(\mathbb{R})$.  
The $\mathbb{R}$-points of $\mathcal{G}$, which is a real Lie group with finitely many connected components, is obtained as:
\[\mathcal{G}(\mathbb{R})= \prod_{\sigma\in S_\infty} G^\sigma(\mathbb{R})\times
\prod_{\{\sigma,\bar{\sigma}\}, \sigma\ne \bar{\sigma}}G^\sigma(\mathbb{C})\cong  (\prod_{\sigma\in S_\infty} G^\sigma(\mathbb{R}))\times G(\mathbb{C})^{d_2} \eqno(3)\] 
where $d_2$ is the number of pairs of conjugate complex embeddings of $F$. 
The second factor is understood to be trivial if $F$ is totally real, that is, if $V_\infty=S_\infty$.

Suppose that $A\in M_n(\mathbb{D})$ is $\tau_r$-hermitian and non-singular.  Assume that $\alpha>0$, $\beta<0$.   
We have an embedding  
$\psi_\mathbb{R}: M_n(\mathbb{D})\rightarrow M_{2n}(\mathbb{R})$, after fixing positive square roots of $\alpha, -\beta$.   See \S \ref{f-structures}. 
Then $\psi_\mathbb{R}(A)$ is symmetric. 
If $\psi_\mathbb{R}(A)$ has exactly $p$ positive eigenvalues, then $\SO(n,\mathbb{H}_\mathbb{R}^{\alpha,\beta})\cong \SO(p, 2n-p)$.   In particular, when $p=0, 2n$, then $A$ is 
definite and the group $\SO(I_n,\mathbb{H}_\mathbb{R}^{\alpha,\beta})\cong \SO(2n)$ is the compact form of $\SO(2n,\mathbb{C})$.   
If $\alpha,\beta <0$, then $\mathbb{H}_\mathbb{R}^{\alpha,\beta}\cong \mathbb{H}$ and $\SO(A,\mathbb{H}^{\alpha,\beta}_\mathbb{R})\cong \SO(n,\mathbb{H})$.

Now suppose that $F$ is a totally real number field, $F\ne \mathbb{Q}$; thus $\sigma(F)\subset \mathbb{R}$ for all $\sigma\in V_\infty.$  
Denote by $\mathbb{D}_\sigma$ the quaternion algebra $\mathbb{H}_{\sigma(F)}^{\sigma(\alpha),\sigma(\beta)}$.   
Suppose that  
$\alpha,\beta\in F$ are such that $\alpha, \beta<0$ and $\sigma(\alpha)>0,\sigma(\beta)<0$ for all $\sigma\in V_\infty, \sigma\ne \iota$.  In view of the fact that $\sigma$ fixes all the rationals, it is easily seen that such elements $\alpha,\beta\in F$ do exist. 

Fix square roots $a,b\in \bar{\mathbb{Q}}$ for $\alpha, -\beta$ respectively 
so that we have an embedding $\psi_{\bar{\mathbb{Q}}}:M_n(\mathbb{D})\to M_{2n}(\bar{\mathbb{Q}})$.   For any 
$\sigma\in V_\infty, \sigma\ne \iota,$ choose an automorphism $\tilde{\sigma}$ of $\bar{\mathbb{Q}}$ 
that extends the isomorphism $\sigma:F\to \sigma(F)$.  We shall denote by $\psi_{\bar{\mathbb{Q}}}^\sigma:M_n(\mathbb{D}_\sigma)\to M_{2n}(\bar{\mathbb{Q}})$ the embedding corresponding to the choice 
$\tilde{\sigma}(a),\tilde{\sigma}(b)$ of square roots of $\sigma(\alpha), \sigma(-\beta)\in \sigma(F)$ respectively.  

Let $A\in M_n(\mathbb{D})$ be a non-singular $\tau_r$-hermitian matrix, that is,  ${}^t\tau_r(A)=A$, $\det(A)\ne 0$.  Since $\tilde{\sigma}(a),\tilde{\sigma}(b)$ are real and since $\sigma(F)\subset \mathbb{R}$, 
the matrix $\psi^\sigma_{\bar{\mathbb{Q}}}(\sigma(A))$ is real symmetric for all $\sigma \ne \iota$.  
Let $\lambda_1,\ldots, \lambda_{2n}\in\bar{ \mathbb{Q}}$ be the eigenvalues of $\psi_{\bar{\mathbb{Q}}}(A)\in M_{2n}(\bar{\mathbb{Q}})$.   Then $\tilde{\sigma}(\lambda_1),\ldots,\tilde{\sigma}(\lambda_{2n})$ are the eigenvalues of $\psi_{\bar{\mathbb{Q}}}^\sigma(\sigma(A))$ for all $\sigma\ne \iota$ in $V_\infty$.  
We choose $A$ such that for {\it all} $\sigma\ne \iota$
the eigenvalues $\tilde{\sigma}(\lambda_j)$ are positive for $1\le j\le 2n$, so that $\psi_{\bar{\mathbb{Q}}}^\sigma(\sigma(A))$ is positive definite. 
In view of the fact that $\tilde{\sigma}(\lambda_j)$ 
is a conjugate of $\lambda_j, 1\le j\le 2n$, it is clear that 
such matrices $A$ exist.   In fact, starting with any $\tau_r$-hermitian matrix $A$, the matrix 
$rI+A\in M_n(\mathbb{D})$ has the required property for any rational number $r>0$ which is larger than the absolute value of any conjugate of 
$\lambda_j, 1\le j\le 2n$.

\begin{theorem} \label{arithmeticlattice}
Let $\iota: F\hookrightarrow  \mathbb{R}$ be a totally real number field, $F\ne \mathbb{Q}$.     Let $\alpha, \beta\in F^\times $.   
Suppose that $\alpha,\beta$ are negative and that $\sigma(\alpha)>0,  \sigma(\beta)<0$ for all $\sigma\ne \iota$ in $V_\infty$.  Let $A\in M_n(\mathbb{H}_F^{\alpha,\beta})$ be $\tau_r$-hermitian and non-singular.  
Suppose that $\psi_{\bar{\mathbb{Q}}}^\sigma(\sigma(A))\in M_{2n}(\mathbb{R})$ is (positive or negative) definite for all $\sigma\ne \iota$.  Then the $\mathcal{O}_F$-points of 
$\SO(A, \mathbb{H}_\mathbb{R}^{\alpha,\beta})\cong \SO(n,\mathbb{H})$ 
is a uniform arithmetic lattice.\\
\end{theorem} 
\begin{proof}  Since $F\ne \mathbb{Q}$, $V_\infty$ contains at least two elements.  
Let $G$ be the $F$-algebraic group $\SO(A,\mathbb{D})=\SO(A, \mathbb{H}^{\alpha,\beta}_F)$.  Our hypotheses that, for $\sigma\in V_\infty, \sigma\ne \iota$,  
$\sigma(\alpha)>0,\sigma(\beta)<0$ and $\psi^\sigma_{\bar{\mathbb{Q}}}(\sigma(A))\in M_{2n}(\mathbb{R})$ is definite, imply that $G^\sigma(\mathbb{R})=\SO(\sigma(A), \mathbb{H}^{\sigma(\alpha),\sigma(\beta)}_{\mathbb{R}})$ is compact. 
Also since $F$ is totally real, it is clear from (3) that 
the only non-compact factor in $\mathcal{G}(\mathbb{R})$ is 
$G(\mathbb{R})=G^\iota(\mathbb{R})=\SO(A,\mathbb{H}^{\alpha,\beta}_\mathbb{R})\cong \SO(n,\mathbb{H}).$  
By the above discussion, we see that $\Gamma:=G(\mathcal{O}_F)$, which is the image of $\mathcal{G}(\mathbb{Z})$ under the projection $\mathcal{G}(\mathbb{R})\to G(\mathbb{R})$ is a lattice in $G(\mathbb{R})$.  
Since  $\mathcal{G}(\mathbb{R})$ has a non-trivial compact factor, $\Gamma$ has no unipotent elements 
and we conclude that $\Gamma$ is uniform. (See \cite{raghunathan}.)
\end{proof}

\begin{remark}{\em 
When $F=
\mathbb{Q}$ the above proof  fails.  
 In fact, let $A=\diag(a_1,\ldots, a_n)$ where the $a_k$ are all non-zero integers.   The $\tau_r$-hermitian inner product on $\mathbb{D}^n$ defined by $A $ when restricted to the subspace $\mathbb{Q}^n$ corresponds 
to the integral quadratic form $B(u):=\langle u,u\rangle=\sum _{1\le k\le n} a_ku_k^2$.  
If not  all the $a_k$ are of the same sign, 
then
 $B$ is an indefinite form and hence it represents a real zero.
It is a classical result that,  if  $n\ge 5$, then $B$ has an {\it integral} zero, say, $u$.  (See \cite[Ch. IV,\S 3]{serre}.)  From this it is not difficult to see that 
$G(\mathbb{Z})$ has unipotent element and hence cannot be a uniform lattice in $G(\mathbb{R})$.  If $a_k>0$ for all $k$, then 
the quadratic form on the $\mathbb{Q}$-vector subspace of $\mathbb{D}^n$ spanned by $ie_1,e_2,\ldots,e_n$ is integral and 
indefinite. Again, as before, it has an integral zero, leading to the conclusion that $G(\mathbb{Z})$ is not uniform.  

We do not know any example of a $\mathbb{Q}$-algebraic group $\mathbf{G}$ with  $G(\mathbb{R})=SO(n,\mathbb{H})$ for which $G(\mathbb{Z})$ is a uniform lattice.
}
\end{remark}

\subsection{Special cycles in $X_\Gamma$}\label{specialcycle}
We now put together the results of \S \ref{restrictionofscalars} and \S \ref{fixedpoints} to construct special cycles 
in $\Gamma\backslash SO(n,\mathbb{H})/U(n)$ for lattices $\Gamma$.  

 Let $A\in M_n(\mathbb{H}_\mathbb{F}^{\alpha,\beta})$ 
be a non-singular $\tau_r$-hermitian matrix where $F\ne \mathbb{Q}$ is any totally real number field and $\alpha,\beta\in F$ satisfy 
the hypotheses of Theorem \ref{arithmeticlattice}.  Thus $\sigma(\alpha)>0,\sigma(\beta)<0$ for all $\sigma\in V_\infty=S_\infty$ except the inclusion $\iota:F\hookrightarrow \mathbb{R}$.  Let $\mathbf{G}=SO(A, \mathbb{H}_F^{\alpha,\beta})$ so that 
$G(\mathbb{R})\cong SO(n,\mathbb{H})$.  With notation as in \S \ref{restrictionofscalars},  
by Theorem \ref{arithmeticlattice}, $G(\mathcal{O}_F)\subset G(\mathbb{R})$ is an arithmetic lattice.  
Our hypotheses on $\alpha,\beta$ and on $A$ imply that $G(\mathcal{O}_F)$ is uniform. 
  
We assume that $A$ is a diagonal matrix $\diag(a_1,\ldots,a_n)$ 
such that $N(a_1)\equiv N(a_l)\mod(F^\times)^2$ for all $l\le n$.  Let $\theta$ denote the Cartan involution of $G(\mathbb{R})$ arising from 
the $F$-Cartan involution  
constructed in Theorem \ref{f-cartan} 
and let $\sigma$ be a sign involution or an involution of even or odd type that commutes 
with $\theta$. (See \S \ref{fixedpoints}.)     
We let $\Gamma$ be any {\it torsionless} finite index subgroup of $G(\mathcal{O}_F)$ that is stable by $\sigma$ and $\theta$.   
In fact any torsionless finite index subgroup of $G(\mathcal{O}_F)$ contains 
such a lattice $\Gamma$.   Let $\Gamma(\sigma)=\Gamma\cap G(\sigma)(\mathbb{R})$.  Clearly it is a discrete subgroup 
of $G(\sigma)(\mathbb{R})$ and we have a closed embedding $\Gamma(\sigma)\backslash G(\sigma)(\mathbb{R})\hookrightarrow 
\Gamma\backslash G(\mathbb{R})$.  It follows that $\Gamma(\sigma)\backslash G(\sigma)(\mathbb{R})$ is compact and so $\Gamma(\sigma)$ is a uniform lattice in $G(\sigma)(\mathbb{R})$.

Let $C(\sigma)=\Gamma(\sigma)\backslash X(\sigma)$ where $X(\sigma)=G(\sigma)(\mathbb{R})/K(\sigma)\subset X:=G(\mathbb{R})/K\cong SO(n,\mathbb{H})/U(n)$.  Since $X(\sigma)$ is connected, so is $C(\sigma)$.  
Thus $C(\sigma)$ is a {\it special cycle} in $X_\Gamma$.   In the case when $\sigma$ is a sign involution or an 
involution of even type, both $C(\sigma), C(\sigma\theta)$ are locally hermitian symmetric and the inclusion 
$C(\sigma)\hookrightarrow X_\Gamma$ is holomorphic.  When $\sigma$ is of odd type, the special cycles 
$C(\sigma), C(\sigma\theta)$ are Lagrangian submanifolds of $X_\Gamma$.     As remarked in the introduction, 
$[C(\sigma)]$, the dual cohomology class of $C(\sigma)$ is non-zero, and, is in fact {\it not} in the image of the 
Matsushima homomorphism $H^*(X_u;\mathbb{C})\to H^*(X_\Gamma;\mathbb{C})$ by \cite[Theorem 2.1]{mr}.

Suppose that $\sigma$ is an involution of even type or of sign type.   Then $C(\sigma)$ is a complex submanifold of 
the compact K\"ahler manifold $X_\Gamma$.  If $C(\sigma)$ is of complex codimension $p$ in $X_\Gamma$, then 
$[C(\sigma)]$ is of Hodge type $(p,p)$, that is, $[C(\sigma)]\in H^{p,p}(X_\Gamma;\mathbb{C})$.    Moreover, 
if $\Lambda$ is a torsionless lattice in $G(\mathbb{R})$ that contains $\Gamma$, as the covering projection 
$\pi: X_\Gamma\to X_\Lambda$ 
is holomorphic, the image of $C(\sigma)$---denoted $C_\Lambda(\sigma)$---is an analytic cycle in $X_\Lambda.$   Again, as $X_\Lambda$ is a 
K\"{a}hler manifold, it represents a non-zero dual cohomology class $[C_\Lambda(\sigma)]\in H^{p,p}(X_\Lambda;\mathbb{C})$ 
which is again not in the image of the Matsushima homomorphism 
$H^*(X_u;\mathbb{C})\to H^*(X_\Lambda;\mathbb{C})$. 
 For, otherwise $[C(\sigma)]=\deg(\pi) .\pi^*([C_\Lambda(\sigma)])\in 
H^*(X_\Gamma;\mathbb{C})$ 
would be in the image of the Matsushima homomorphism, in view of the following commutative diagram 
\[
\begin{array}{ccc}
H^*(X_u;\mathbb{C})& \to & H^*(X_\Lambda;\mathbb{C})\\
id \downarrow & & \downarrow  \pi^*\\
H^*(X_u;\mathbb{C}) &\to & H^*(X_\Gamma; \mathbb{C})\\
\end{array}
\]
where the horizontal arrows are the Matsushima homomorphisms. 

Next assume that $\sigma$ is of odd type and $n$ is odd. Then, by Lemma \ref{orforoddtype}, $\sigma$ satisfies 
the condition {\em Or} of \cite[Theorem 4.11]{rs}.   Therefore there exists a finite index subgroup $\Lambda\subset \Gamma$ which is stable by $\sigma$ and $\theta$ such that the submanifold $C_\Lambda(\sigma)$, more briefly denoted $C(\sigma)$, defined as  $C(\sigma):=\Lambda\backslash G(\sigma)(\mathbb{R})/K(\sigma)$ is a special cycle and $[C(\sigma)] \cdot [C(\sigma\theta)]\ne 0$ in $H^*(X_\Lambda;\mathbb{C})$.  
Note that $\dim C(\sigma)={n\choose 2}=(1/2)\dim X_\Gamma$ so that $[C(\sigma)]\in H^{{n\choose 2}}(X_\Lambda;\mathbb{C})$.  

\begin{example} \label{extraspecialcycle}
Taking $\sigma$ to be the involution corresponding to $S_{1,n-1}=\diag(1,-1,\ldots,-1)$, we obtain that $G(\sigma)(\mathbb{R})\cong SO(1, \mathbb{H})\times SO(n-1,\mathbb{H})\cong \mathbb{S}^1\times SO(n-1,\mathbb{H}) $ and the special cycle $C_\Lambda(\sigma)$ has (complex) codimension $n-1$ in $X_\Lambda$.  
\end{example}

\section{$L^2(\Gamma\backslash G)$ and $\theta$-stable parabolic subalgebras} \label{thetastable}

\subsection{Relative Lie algebra cohomology and Matsushima isomorphism}\label{gkcohomology}
Let $\Gamma$ be a uniform lattice in a linear connected semisimple group $G$.  We keep the notations of \S \ref{intro}.  Thus $K$ denotes  
a maximal compact subgroup of $G$, $\theta$ the corresponding Cartan involution, and $\mathfrak{g}=\mathfrak{k}\oplus \mathfrak{p}$ the Cartan decomposition.  

Recall that a $\theta$-stable parabolic subalgebra of $\mathfrak{g}_0$ 
is a parabolic subalgebra $\mathfrak{q}$ of $\frak{g}$ such that $\theta(\mathfrak{q})=\frak{q}$ and 
$\bar{\mathfrak{q}}\cap \mathfrak{q}=:\mathfrak{l}$ is the  Levi subalgebra of $\mathfrak{q}$.   Here, the bar refers to conjugation of $\mathfrak{g}=\mathfrak{g}_0\oplus i\mathfrak{g}_0$ with respect to the real form $\mathfrak{g}_0$. 
As recalled in \S\ref{intro}, a result of Gelfand and Piatetskii-Shapiro \cite{gp}, \cite{ggp} states that $L^2(\Gamma\backslash G)$ 
decomposes into a Hilbert sum $\hat{\oplus}_{\pi\in \hat{G}} m(\pi, \Gamma) H_\pi$ of irreducible 
unitary representations $H_\pi$ with finite multiplicities $m( \pi,\Gamma)$.   One has the Matsushima isomorphism $H^*(\Gamma;\mathbb{C})\cong H^*(\mathfrak{g},K;L^2(\Gamma\backslash G)_K)=\oplus_{\pi\in \hat{G}}m(\pi,\Gamma) H^*(\mathfrak{g},K;H_{\pi,K})$. 
Since $ \Gamma\backslash G$ is compact, there are only finitely many $\pi\in \hat{G}$ with $m(\pi,\Gamma)\ne 0$ having non-vanishing $(\mathfrak{g},K)$-cohomology. In fact, to each  
$\theta$-stable parabolic subalgebra $\mathfrak{q}$ of $\mathfrak{g}_0$ there is an irreducible unitary $G$-representation $(\mathcal{A}_{\mathfrak{q}},A_\mathfrak{q})$ associated to $\mathfrak{q}$ with $H^*(\mathfrak{g},K;A_{\mathfrak{q},K})\ne 0$. If 
$V$ is a $(\mathfrak{g},K)$-module with $H^*(\mathfrak{g},K;V)\ne 0$, then $V$ is isomorphic, as a $(\mathfrak{g},K)$-module to the space of (smooth) $K$-finite vectors of $A_\mathfrak{q}$ for some $\theta$-stable parabolic subgalgebra of $\mathfrak{g}_0$. 
Up to unitary equivalence, there are only 
finitely many representations of the form $(\mathcal{A}_\mathfrak{q},A_\mathfrak{q})$.  Hence   
we have 
\[ H^*(\Gamma; \mathbb{C})\cong \bigoplus_{[\mathfrak{q}]}m(\mathcal{A}_\mathfrak{q},\Gamma)H^*(\mathfrak{g}, K;A_{\mathfrak{q},K})\]
where the sum is over all equivalence classes of $\theta$-stable parabolic subalgebras $\mathfrak{q}$ of $\mathfrak{g}_0$;  
here $\mathfrak{q}$ and $\mathfrak{q}'$ belong to the same equivalence class $[\mathfrak{q}]$ if $\mathcal{A}_{\mathfrak{q}}$ and $\mathcal{A}_{\mathfrak{q}'}$ are unitarily equivalent.   

It is known that two irreducible unitary representations $(\pi, H_\pi), (\sigma, H_\sigma)$ of $G$ 
are unitarily equivalent if their Harish-Chandra modules $H_{\pi,K}$ and $H_{\sigma,K}$ are isomorphic as $(\mathfrak{g},K)$-representations.  Moreover, any irreducible unitary $(\mathfrak{g},K)$-representation arises as the Harish-Chandra module of a {\it 
unique} irreducible unitary $G$-representation.  (See \cite[Ch. IX]{knapp}.)
 Thus, in order to describe $\mathcal{A}_{\mathfrak{q}}$, it suffices to describe 
its Harish-Chandra module.  We shall do this in the special case when the (complex) rank of $G$ equals the rank of $K$ as this condition 
holds in the case when $G=SO(n,\mathbb{H})$.   Fix a maximal 
torus $T\subset K$. In view of our assumption on $G$, $\mathfrak{t}=\text{Lie}(T)\otimes \mathbb{C}=\mathfrak{t}_0\otimes{C}$ is a Cartan subgalgebra of $\mathfrak{g}$.  We assume that 
$\mathfrak{t}\subset \mathfrak{q}$.  Write $\mathfrak{q}=\mathfrak{l}\oplus\mathfrak{u}$ where $\mathfrak{u}$ is the nilradical of $\mathfrak{q}$. 
Choose a positive system $(\mathfrak{l}\cap \mathfrak{k},\mathfrak{t})$ and extend it to a positive system 
for $(\mathfrak{k},\mathfrak{t})$ so that the weights of $\mathfrak{u}\cap \mathfrak{k}$ are 
all positive roots of $\mathfrak{k}$.  
Then $\mathcal{A}_{\mathfrak{q}}$ is determined, up to unitary equivalence, by the set of weights in $\mathfrak{u}\cap\mathfrak{p}$.  (See Remark{equivalentq}(i).)

Now let $L\subset G$ be the Lie subgroup corresponding to the Lie subalgebra $\mathfrak{l}_0:=\mathfrak{l}\cap \mathfrak{g}_0.$ Then $K\cap L$ is a maximal compact subgroup of $L$.  Let $Y_\mathfrak{q}$ denote the compact dual of $L/K\cap L$.  It turns out that $H^r(\mathfrak{g},K;A_{\frak{q},K})=\hom_K(\Lambda^r(\mathfrak{p}),A_{\mathfrak{q},K}) \cong H^{r-R(\frak{q})}(Y_\mathfrak{q};\mathbb{C})$, where $R(\mathfrak{q}):=\dim_\mathbb{C}(\mathfrak{p}\cap \mathfrak{u})$.  

When $rank(K)=rank(G)$ and $\mathfrak{q}$ is a Borel subalgebra,  $\mathcal{A}_\mathfrak{q}$ is a discrete series 
representation.  In this case, $L=T$ and $R(\mathfrak{q})=(1/2) \dim G/K$ and $Y_\mathfrak{q}$ is a point. 

Suppose that $G/K$ is a Hermitian symmetric space, equivalently the centre of $T$ is non-discrete. Then the space $X_\Gamma=\Gamma\backslash G/K$ admits the structure of a smooth projective variety arising from a $G$-invariant complex structure on $G/K$.  The tangent space $\mathfrak{p}_0$ being a complex vector space, $\mathfrak{p}=\mathfrak{p}_+\oplus \mathfrak{p}_-$ where $\mathfrak{p}_+$ and $\mathfrak{p}_-$ are conjugate complex vector spaces.
 The real tangent space $\mathfrak{p}_0$ may be identified with the holomorphic tangent space 
$\mathfrak{p}_+$ as $K$-representations.  Note that $\mathfrak{p}_-$ is the dual of $\mathfrak{p}_+$ as a $K$-module.  Also 
$\mathfrak{p}_+,\mathfrak{p}_-$ are abelian subalgebras of $\mathfrak{g}$.
The Hodge structure on $H^r(\Gamma;\mathbb{C})$ arises from the decomposition $\Lambda^r(\mathfrak{p})
=\Lambda^r(\mathfrak{p}_+\oplus \mathfrak{p}_-) 
=\oplus_{a+b=r}\Lambda^a(\mathfrak{p}_+)\otimes \Lambda^b(\mathfrak{p}_-)$.    More precisely,  $H^{a,b}(\mathfrak{g},K;A_{\mathfrak{q},K})\cong \hom_K(\Lambda^a(\mathfrak{p}_+)\otimes \Lambda^b(\mathfrak{p}_-),A_{\mathfrak{q},K})$ and we have 
\[H^{a,b}(\Gamma;\mathbb{C})=\oplus_{a+b=r} m(\mathfrak{q},\Gamma) H^{a,b}(\mathfrak{g},K; 
A_{\mathfrak{q},K}), \eqno(4)\]
where $m(\mathfrak{q},\Lambda)$ stands for $m(\mathcal{A}_\mathfrak{q},\Lambda)$.  
  See \cite[Ch. VII, \S\S2,3]{bw}.
Let $R_\pm(\mathfrak{q})=\dim_\mathbb{C}(\mathfrak{p}_\pm\cap \mathfrak{u})$ so that $R(\mathfrak{q})=R_+(\mathfrak{q})+R_-(\mathfrak{q})$.  Then 
\[H^{a,b}(\mathfrak{g},K;A_{\mathfrak{q},K})=0~\textrm{if} ~a-b\ne R_+(\mathfrak{q})-R_-(\mathfrak{q}).\eqno(5)\]  See \cite[Proposition 6.19]{vz}.    

The construction of the Harish-Chandra module of $\mathcal{A}_{\mathfrak{q}}$ was originally due to Parthasarathy \cite{parthasarathy}.    
Vogan-Zuckerman \cite{vz}, Vogan \cite{vogan} gave a construction in terms of cohomological induction and showed that they are unitarizable.  A very readable 
account explaining the basic properties of $\mathcal{A}_\mathfrak{q}$ is given in \cite{vogan97}. 

\subsection{The $\theta$-stable parabolic subalgebras of $\SO(n,\mathbb{ H})$.}\label{thetastableparabolic}

Recall from  \S\ref{realform} that $\SO(n,\mathbb{H})$
is a real form of $\SO(2n,\mathbb{C})$.  We shall identify $\SO(n,\mathbb{H})$ with its image under $\psi: \SO(n,\mathbb{H})\to \SO(2n,\mathbb{C})$.  Write $G:=\SO(n,\mathbb{H}), K=\{\left ( \begin{matrix}Z & W \\-W& Z\end{matrix}\right)\in G\mid Z,W\in M_n(\mathbb{R}) \}\cong U(n)$, a maximal compact subgroup of $G$. (See \S \ref{realform}.)  Let $T\subset K$ 
be the maximal torus $\{\left(\begin{matrix}Z & W \\-W & Z\end{matrix}\right)\in K\mid Z, W ~\textrm{diagonal}\}$.  
Note that $rank(K)=rank(G(\mathbb{C}))$.  So $\mathfrak{t}$ is a Cartan subalgebra of $\frak{g}=\frak{so}(2n,\mathbb{C})$. 

Denote by $\underline{B}$ the matrix $\left(\begin{matrix}0 & B\\-B & 0\end{matrix}\right)$ where $B\in M_n(\mathbb{C}).$ 
The Lie algebra $\mathfrak{t}_0$ of $T$ is 
$\{\underline{B}\in M_{2n}(\mathbb{R})\mid B\textrm{~is diagonal}.\}$.   
We set $\epsilon_j:\frak{t}\to \mathbb{C}$ be the $\mathbb{C}$-linear form defined as 
$\epsilon_j(\underline{B})=-ib_j\in \mathbb{C}$ where $B=\diag(b_1,\ldots,b_n)$.
Note that $\epsilon_j$ takes real values on $i\mathfrak{t}_0$. 
 Then 
$\Phi=\{\pm(\epsilon_i\pm\epsilon_j)\mid 1\le i<j\le n\}, \Phi_\mathfrak{k}=\{\epsilon_i-\epsilon_j\mid 1\le i, j\le n, ~i\ne j\}$ 
are the set of roots of $\mathfrak{g}$, respectively, of $\mathfrak{k}$.   
We take $\Phi^+$ to be $\{\epsilon_i\pm\epsilon_j\mid 1\le i<j\le n\}$ and $\Phi^+_\mathfrak{k}=\Phi_\mathfrak{k}\cap 
\Phi^+.$   We denote by $\Phi_n$ the set of non-compact roots  $\{\pm(\epsilon_p+\epsilon_q)\mid 1\le p<q\le n\}$ and by $\Phi^+_n$ the 
set $\Phi^+\cap \Phi_n$ of positive non-compact roots.

Let $x\in i\mathfrak{t}_0,$ where $x:=i\underline{X}, X=\diag(x_1,\ldots,x_n)\in M_n(\mathbb{R})$.  
Applying $Ad(g)$ for a suitable $g\in K$ if necessary, we may (and do) assume that 
$x_1\ge\ldots\ge x_n$, so that 
$\alpha(x)\ge 0~\forall\alpha\in \Phi^+_\frak{k}$.     
If $(\epsilon_p+\epsilon_q)(x)\ge 0, p>q$, then 
$(\epsilon_r+\epsilon_s)(x)\ge 0~\forall r\le p, s\le q.$  
The centralizer $\mathfrak{l}_x$ of 
$x\in\mathfrak{g}$ is a reductive subalgebra of $\mathfrak{g}$ that contains $\mathfrak{t}$.  
It contains a root space $\mathfrak{g}_\alpha$ whenever (a) $\alpha=\epsilon_i-\epsilon_j$ and $x_i=x_j$, or 
(b) $\alpha=\pm(\epsilon_i+\epsilon_j)$ and $x_i=-x_j$.   
 
We have $\mathfrak{g}=\mathfrak{l}_x\oplus \mathfrak{u}\oplus 
\mathfrak{u}_-$, where $\mathfrak{u}=\mathfrak{u}_x$ is the nilpotent subalgebra 
$\oplus_{\alpha\in \Phi, \alpha(x)>0}\mathfrak{g}_\alpha$ 
and $\mathfrak{u}_-=\oplus_{\alpha\in \Phi,\alpha(x)<0} \mathfrak{g}_{\alpha}$.    
Then $\mathfrak{q}=\mathfrak{q}_x:=\mathfrak{l}_x\oplus \mathfrak{u}_x$ is a 
$\theta$-stable parabolic subalgebra of $\mathfrak{g}_0$.  Up to conjugation by 
$K$, all $\theta$-stable parabolic subalgebras in $\mathfrak{g}$ arise this way.

Consider the equivalence relation $p\sim q$ if $|x_p|=|x_q|$ 
on the set $\{1\le j\le n\}$.  If $x_p=0$ we denote the corresponding equivalence class $[p]$ by $\mathcal{N}_x$ and its cardinality $\#\mathcal{N}_x$ by $N_x$. If 
$x_p\ne 0$, we define two subsets $I_{[p]}:=\{ j\in [p]\mid  x_j> 0\}, 
J_{[p]}:=\{j\in [p]\mid  x_j<0\}$ of $[p]$.  Note that $I_{[p]}, J_{[p]}$ 
are disjoint sets of consecutive integers, at least one of which is non-empty and $I_{[p]}\cup J_{[p]}=[p]$ .   
The sets $\mathcal{N}_x, I_{[p]}, J_{[p]}, 1\le p\le n, x_p\ne 0$, form a partition of the integers $1$ up to $n$. 
Denote by $L_x$ the set $\{\max I_{[p]}\mid\exists q~ x_p=-x_q\},$ and by $R_x$ the set $\{\min J_{[p]}\mid \exists q ~x_p=-x_q\}$.    
Clearly there is a bijection $L_x\to R_x$ sending $p\in L_x$ to the unique $p'\in R_x$ such that 
$x_p=-x_{p'}$. 
Also set $m:=\max \mathcal{N}_x$ if $\mathcal{N}_x$ is non-empty. 
We shall now describe the root system of $\mathfrak{l}$ with respect to $\mathfrak{t}$.  
Let $\Phi_x:=\{\epsilon_i-\epsilon_j\mid x_i=x_j\}\cup 
\{\pm(\epsilon_i+\epsilon_j)\mid x_i=-x_j\}$ and let $\Delta_x= 
\Delta_x^+\cup \Delta^-_x\cup \Delta^0$ where 
$\Delta^+_x:=\{\epsilon_i-\epsilon_{i+1} \mid x_i=x_{i+1}\ge 0\}
\cup \{\epsilon_i-\epsilon_{i+1}\mid i,i+1\in J_{[i]}, ~I_{[i]}=\emptyset \},  
\Delta^-_x:=\{\epsilon_{i+1}-\epsilon_i\mid i,i+1\in J_{[i]}, ~I_{[i]}\ne \emptyset\}$ and $\Delta^0_x:=\{\epsilon_p+\epsilon_{p'}\mid p\in L_x\} \cup \{ \epsilon_{m-1}+\epsilon_m\mid x_{m-1}=x_m= 0\}$; thus  
$\epsilon_{m-1}+\epsilon_m\in \Delta^0_x$ if and only if  $N_x\ge 2$. 
Then $\mathfrak{l}_x=\mathfrak{t}\oplus _{\alpha\in \Phi_x}\mathfrak{g}_\alpha$.   Also, $\Delta_x$ is a set of 
simple roots for a positive 
root system $\Phi^+_x\subset \Phi_x$.  
Write $\Phi_{[p]}=\{\alpha\in \Phi_x\mid \alpha=\pm (\epsilon_i\pm \epsilon_j),~
i,j\in [p]\}$ and  $\Delta_{[p]}= \Delta_x\cap \Phi_{[p]}$.  Then $(\Phi_{[p]},\Delta_{[p]})$ is an irreducible root system except when 
$x_p=0, $ and $\#[p]=N_x=2$.  Also
$\Phi_x=\bigcup_{[p]}\Phi_{[p]}$ with $\Phi_{[p]}, \Phi_{[q]}$ being orthogonal if $[p]\ne [q].$   
Moreover $[\mathfrak{l}_x,\mathfrak{l}_x]$ is a direct sum of simple ideals $\mathfrak{s}_{[p]}$ whose root system 
(relative to $\mathfrak{t}_{[p]}:=\mathfrak{t}\cap \mathfrak{s}_{[p]}$) is given by restriction of elements of 
$\Phi_{[p]}$ to $\mathfrak{t}_{[p]}$.  The Killing-Cartan type of the 
Lie algebra $\mathfrak{s}_{[p]}$ is (a) $A_{\#[p]-1}$ if $x_p\ne 0$, (b)  
 type $D_{N_x}$ when $N_x\ge 2$. (Of course, a genuine type $D$ factor occurs only when $N_x\ge 4$.)
   The radical of $\mathfrak{l}_x$ is isomorphic to $\mathbb{C}^s$ where $s$ equals the number of singletons among the sets  $\mathcal{N}_x, [p], 1\le p\le n$.

Let $R(\mathfrak{q})=\dim_\mathbb{C}(\mathfrak{p}\cap \mathfrak{u})$.  As remarked in \S\ref{gkcohomology}, 
the first non-vanishing $(\mathfrak{g},K)$-cohomology group $H^j(\mathfrak{g},K;A_{\mathfrak{q},K})$ 
occurs in dimension $j=R(\mathfrak{q})$.  
The following proposition gives a combinatorial formula for $R(\mathfrak{q})$.

\begin{proposition}\label{idealsofl}
Let $x=(x_1,\ldots, x_n)$ where $x_i\ge x_j, 1\le i<j\le n$.  With the above notations, we have 
\[[\mathfrak{l}_x,\mathfrak{l}_x]\cong  \oplus_{[p],x_p\ne 0} \mathfrak{sl}(\#[p],\mathbb{C})\oplus \mathfrak{so}(2N_x,\mathbb{C})\]
where the last summand occurs only when $N_x\ge 2$.  
Also, 
\[R(\mathfrak{q}_x)=\dim_\mathbb{C}(\mathfrak{u}\cap\mathfrak{p})={n\choose 2}-{N_x\choose 2}- \sum_{[p], x_p\ne 0}\#I_{[p]}\cdot\#J_{[p]}.\]  
\end{proposition}
\begin{proof}
In view of the above discussion, we need only establish the asserted formula for $R(\mathfrak{q}_x)$.
Note that $R(\mathfrak{q}_x)
=\#\Phi^+_n-\#\{\alpha\in \Phi^+_n\mid \alpha, -\alpha\in \Phi_x\}$.   Note that $\alpha,-\alpha\in \Phi_x$ if and only 
if $\alpha\in \Phi(\mathfrak{l},\mathfrak{t})$.  
So we need to count the number of $\alpha\in \Phi_n^+$ 
which are roots of $\mathfrak{l}$.  
Let  $\alpha:=\epsilon_i+\epsilon_j$ and let $\alpha,-\alpha \in \Phi_x, i<j$; thus $x_i=-x_j$. 

{\it Case 1.} If $x_i=0=x_j$, the root space $\mathfrak{g}_\alpha
$ is contained in the ideal corresponding to the summand $\mathfrak{so}(2N_x,\mathbb{C})$.  There are exactly 
${N_x\choose 2}$ many such positive roots.  

{Case 2.}  Let  $x_i\ne 0$.  Then $0=\alpha(x)=x_i+x_j=0$ implies that $x_i>0, x_j<0$ and so there exists a unique $p\in L_x$ such that $p\sim i\sim j$. So $k\in I_{[p]}$ for $i\le k\le p$ and $l\in J_{[p]}$ 
for $p'\le l\le j$.   
We see that $\alpha=\sum_{i\le k<p} (\epsilon_k-\epsilon_{k+1})+\epsilon_p+\epsilon_{p'}+\sum_{p'\le l< j}
(\epsilon_{l+1}-\epsilon_l)$ is an expression of $\alpha$ as a sum of simple roots of $(\Phi_x,\Delta_x)$.  Hence $\alpha\in \Phi^+_x$ and it 
occurs in the simple ideal $\mathfrak{s}_{[p]}$ isomorphic to $\mathfrak{sl}(\#I_{[p]}+\#J_{[p]})=\mathfrak{sl}(\#[p])$.  
It is clear that there are exactly $\#I_{[p]}\cdot \#J_{[p]}$ such positive 
roots in $\mathfrak{sl}([p])$.  
It follows that the set $\{\alpha \in \Phi^+_n\mid \alpha\in \Phi(\mathfrak{l},\mathfrak{t})\}$ has cardinality equal to ${N_x\choose 2}+\sum_{[p], x_p\ne 0} \#I_p\cdot \#J_p$.  This completes the proof.
\end{proof}

Recall that $R_\pm(\mathfrak{q})=\dim_\mathbb{C}(\mathfrak{u}\cap \mathfrak{p}_\pm)$.   The following 
lemma provides a formula for $R_+(\mathfrak{q}_x)$.  An analogous formula for $R_-(\mathfrak{q})$ can be 
readily obtained.  We omit the proof of this lemma, which is immediate from the definitions 
of $I_{[p]}, J_{[p]}$.

\begin{lemma}
Let $x=(x_1,\ldots, x_n)$ be a monotone decreasing sequence of real numbers.   For $1\le p<n, x_p>0$,  
 let $q(p)$ denote the largest integer in $[p,n]$ such that $x_{q(p)}>-x_p$. 
We have $R_+(\mathfrak{q}_x)=\sum_{x_p>0, 1\le p<n} q(p)-p$.   
\hfill $\Box$
\end{lemma}
Suppose that $x_p>0$.   If $J_{[p]}\ne \emptyset$, then $q(p)=\min(J_{[p]})-1$; if $J_{[r]}=\emptyset$ for all $r\ge p$, then $q(p)=n$.

\begin{lemma} \label{duality}
Let $x=(x_1,\ldots, x_n)$ be a monotone decreasing sequence of real numbers.  Let $\iota(x):=(-x_n,\ldots, -x_1).$   Then 
$R_\pm(\mathfrak{q}_x)=R_\mp(\mathfrak{q}_{\iota(x)})$.  
\end{lemma}
\begin{proof}
It is readily checked that the root space $\mathfrak{g}_{\iota(\alpha)}\subset \mathfrak{q}_{\iota(x)}$ if and only if $\mathfrak{g}_{-\alpha} \subset 
\mathfrak{q}_x$ where $\iota:\Phi\to \Phi$ is the bijection induced by $\iota(\epsilon_p)=\epsilon_{n+1-p}$.  In particular, when $i\ne j$, 
$\mathfrak{g}_{\epsilon_{n+1-j}+\epsilon_{n+1-i}}\subset \mathfrak{u}_{\iota(x)}$ if and only if $-x_j-x_i>0$ if and only if $x_i+x_j<0$ if and only 
if $\mathfrak{g}_{-\epsilon_i-\epsilon_j}\subset \mathfrak{u}_x$.  Since $\mathfrak{g}_{\iota(\alpha)}\subset \mathfrak{p}_+$ if and 
only if $\mathfrak{g}_{-\alpha}\subset \mathfrak{p}_-$, it follows that $R_+(\mathfrak{q}_{\iota(x)})=R_-(\mathfrak{q}_x)$. This completes the 
proof.
\end{proof}

In the examples below, we often write $R_+$ for $R_+(\mathfrak{q})$, etc.,  when there is 
no possibility of confusion.  

\begin{example} \label{minimalparabolicseg}
{\em 
(i) Let $s=s(k,l)=(1,\ldots, 1, 0\ldots,0,-1, -1,\ldots, -1)$, where there are $k$ many $1$s and $l$ many $-1$s.  Then 
$R_+(\mathfrak{q}_s)=(n-l-1)+(n-l-2)+\ldots +(n-l-k)=k(n-l-k)+{k\choose 2}$; similarly $R_-=l(n-k-l)+{l\choose 2}.$ \\
(ii) More generally, if $x_1=-x_n>0$, then $R_+\ge k(n-l-k)+{k\choose 2}$ where $k=\# I_{[1]}, l=\#J_{[1]}
\ge 1$. 
}
\end{example}

Recall that $Y_\mathfrak{q}$ is the compact dual of the symmetric space $L/L\cap K$ where $L\subset G$ is the 
Lie subgroup whose Lie algebra is $\mathfrak{l}_0=\mathfrak{l}\cap \mathfrak{g}_0\subset \mathfrak{g}_0$.   
The roots of $\mathfrak{k}_0\cap \mathfrak{l}_0$ with respect to $\mathfrak{t}_0$ are precisely the compact roots in $\Phi_x$ (restricted to $\mathfrak{t}_0$).  The root space of $\mathfrak{k}_0$ corresponding to a compact root $\alpha$ will be 
denoted $\mathfrak{k}_{0,\pm\alpha}:=\mathfrak{k}_0\cap(\mathfrak{k}_\alpha+\mathfrak{k}_{-\alpha})$. 
Thus   $\mathfrak{k}_0\cap \mathfrak{l}_0=\mathfrak{t}_0\oplus(\oplus_{\pm\alpha} \mathfrak{k}_{0,\pm\alpha})$, where the sum is over pairs of roots $\pm\alpha$ in $\Phi_x\cap \Phi_\mathfrak{k}$.
Therefore the set of compact roots of $[\mathfrak{l}_0,\mathfrak{l}_0]$ are the restrictions of elements of 
$\Phi_x\cap \Phi_\mathfrak{k}$ to $\mathfrak{t}_0\cap [\mathfrak{l}_0,\mathfrak{l}_0]$.  We shall make no distinction 
in the notation between a root of $(\mathfrak{k}_0,\mathfrak{t}_0)$ and its restriction to $\mathfrak{t}_0\cap [\mathfrak{l}_0,\mathfrak{l}_0]$.

\begin{lemma} \label{derivedl0}
The derived algebra $[\mathfrak{l}_0,\mathfrak{l}_0]$ is a direct sum of the following simple Lie algebras 
$\mathfrak{s}_{[p],0}$, the sum being 
over the set of all equivalence classes $[p]\subset \{1,2,\ldots, n\}$:\\
(i) $\mathfrak{su}(\#[p])$ if $I_{[p]}$ or $J_{[p]}$ is empty,\\
(ii) $\mathfrak{su}(\#I_{[p]},\#J_{[p]})$ if both $I_{[p]}, J_{[p]}$ are non-empty,\\
(iii) $\mathfrak{so}(2N_x,\mathbb{H})$, when $N_x\ge 2$.\hfill $\Box$
\end{lemma}

We omit the proof.  In fact, it is possible to describe the simple factors of the reductive Lie group $L$ which implies 
the above lemma.  Alternatively an easy computation yields the compact and non-compact roots of the real form $\mathfrak{s}_{[p],0}$ 
of $\mathfrak{s}_{[p]}$, which allows one to determine the Lie algebra $\mathfrak{s}_{[p],0}$.

We are now ready to describe the symmetric space $Y_\mathfrak{q}$, the compact dual of the symmetric space 
$L/L\cap K \cong [L,L]/([L,L]\cap K)$.  We have $Y_\mathfrak{q}=M/[L,L]\cap K$ where $M$ is the maximal compact subgroup 
of the complex semisimple Lie group with Lie algebra $[\mathfrak{l},\mathfrak{l}]$. 
From Proposition \ref{idealsofl} and the above lemma we see that $M=(\prod_{[p], x_p\ne 0} SU(\#[p]))\times \SO(2N_x)$, 
$[L,L]\cap K=(\prod_{[p], x_p\ne 0} S(U(\#I_{[p]})\times U(\#J_{[p]})))\times U(N_x)$.
It is understood that $U(0)$ is trivial and that the last factor in $M$ and $[L,L]\cap K$ is present only if $N_x\ge 2$.
Hence $Y_\mathfrak{q}=(\prod_{[p],x_p\ne 0} \mathbb{C}G_{\#[p], \#I_{[p]}})\times \SO(2N_x)/U(N_x)$, where $\mathbb{C}G_{k+l,l}$ denotes the complex Grassmann manifold of $l$-planes in $\mathbb{C}^{k+l}$. 

In view of the fact that $H^j(\mathfrak{g},K;A_{\mathfrak{q},K})\cong H^{j-R(\mathfrak{q})}(Y_{\mathfrak{q}};\mathbb{C})$ we obtain an explicit formula for the Poincar\'e polynomial of  $H^*(\mathfrak{g},K;A_{\mathfrak{q},K})$.

Recall, from \cite{borel}, that the Poincar\'e polynomial $P_t(X)$ (where cohomology with coefficients in a field of characteristic zero is understood) of a homogeneous manifold of the form $X=M/H$  where $M$ is a compact connected 
Lie group and $H$ is a connected subgroup having the same rank as $M$ is given by the formula of Hirsch.   
This is applicable to 
the complex Grassmann manifold $\mathbb{C}G_{k+l,k}$ and to $\SO(2k)/U(k)$.  We have 
\[P_t(\mathbb{C}G_{k+l,k})=(1-v^{l+1})\cdots (1-v^{k+l})/((1-v)\cdots(1-v^k))\]
where $v:=t^2$.  (See \cite{milnor-stasheff}.)   Also, the Poincar\'e polynomial of $\SO(2k)/U(k)$ is 
\[P_t(\SO(2k)/U(k))=
(1+v)(1+v^2)\cdots (1+v^{k-1}).
\]

Recall that $R=R(\mathfrak{q})=\dim_\mathbb{C}(\mathfrak{u}_x\cap \mathfrak{p})$.

\begin{theorem}  \label{poincarepolynomial}
We keep the notations as above and set $\mathfrak{q}:=\mathfrak{q}_x$. 
The Poincar\'e polynomial $P_t(\mathcal{A}_\mathfrak{q})$  of $H^*(\mathfrak{g}, K; A_{\mathfrak{q},K})$ is given by 
\[P_t(\mathcal{A}_\mathfrak{q})=t^{R(\mathfrak{q})}(\prod_{[p], x_p\ne 0} P_t(\mathbb{C}G_{\#[p], \#I_{[p]}}))\times P_t(\SO(2N_x)/U(|N_x|)).\] \hfill $\Box$
\end{theorem}

\subsection{Representations $\mathcal{A}_\mathfrak{q}$ with $R_\pm(\mathfrak{q})\le n-1$}

We consider the equivalence relation on the set of all $\theta$-stable parabolic subalgebras of $\mathfrak{g}_0=\mathfrak{so}(n,\mathbb{H})$ defined as $\mathfrak{q}\sim \mathfrak{q}'$ if $\mathcal{A}_\mathfrak{q}$ and $\mathcal{A}_{\mathfrak{q}'}$ 
are unitarily equivalent, or equivalently, if the corresponding Harish-Chandra modules are isomorphic as $(\mathfrak{g},K)$-modules.  
Clearly $\mathfrak{q}\sim \mathfrak{q}'$ if they belong to the same $K$-conjugacy class.  Thus each conjugacy 
class is represented by a $\theta$-stable parabolic subalgebras of the form $\mathfrak{q}_x$ with $x=(x_1,\ldots,x_n)\in i\mathfrak{t}_0$ monotone 
decreasing.  Two such elements $x,x'$ define equivalent representations $\mathcal{A}_{\mathfrak{q}_x},\mathcal{A}_{\mathfrak{q}_{x'}}$ if and only if $\mathfrak{u}_x\cap \mathfrak{p}=\mathfrak{u}_{x'}\cap\mathfrak{p}$. (See Remark \ref{equivalentq}(i).)
We shall classify all $\mathfrak{q}_x$ for which $R_+(\mathfrak{q}_x),R_-(\mathfrak{q}_x)\le n-1$ assuming $n>8$.  In view of Lemma \ref{duality} 
we assume without loss of generality that $R_-(\mathfrak{q}_x)\le  R_+(\mathfrak{q}_x)\le n-1$.  
We shall write $R_\pm$ to denote $R_\pm(\mathfrak{q}_x)$.  We ignore the case $R_+=R_-=0$ which corresponds to $\mathfrak{q}=\mathfrak{g}$. So $x_1>0.$

\begin{proposition} \label{minimalparabolics} Let $n>8$.
Let $x=(x_1,\ldots,x_n)$ be a monotone decreasing with $x_1>0$.  Suppose that $0\le R_-\le R_+\le n-1$.  Then the following statements hold:\\
(i) If $R_+>0$, then $R_+\ge n-2$. \\
(ii) If $R_+=n-1$, then $R_-\in\{0,n-2\}$. Moreover $\mathfrak{q}_x=\mathfrak{q}_s$ where $s=(1,0,\ldots,0,0)$ or $(2,0\ldots,0,-1)$.\\
(iii) If $R_+=n-2$, then $R_-=n-2$ and $x_1=-x_n, x_j=0, 2\le j<n$.  
\end{proposition}  
\begin{proof}
The proof is broken up into three cases depending on whether $x_1+x_n$ is positive, negative, or zero. 
 
{\it Case 1:} 
Suppose that $x_1+x_n>0$.  Then $R_+\ge n-1$ with equality if and only if $x_2+x_3\le 0$.  Our assumption that $R_+\le n-1$ implies that we 
must have $R_+=n-1$ and $x_2+x_{n-1}\le 0$.  
If 
$ x_2+x_{n-1}< 0$, then $x_i+x_j<0, ~2\le i<n-1\le j\le n$. In particular $R_-\ge 2n-5$. Since $R_-\le n-1$, this is impossible if  $n>4$.  
So $x_2+x_{n-1}=0$.  This implies that $x_2=-x_j, 3\le j\le n-1$.   

If $x_2\ne 0$, then $x_j<0~\forall j\ge 3$.  Hence $x_i+x_j<0$ for 
$3\le i<j\le n$.  Thus $R_-\ge {{n-2}\choose{2}}\ge 2n-7$. As $R_-\le n-1$, we must have $n-1\ge 2n-7$ contradicting $n>6$. 
So $x_2=0$.  This means that $x_j=0~\forall j\le n-1$.  Now $x_1+x_n>0$ implies that $R_-=n-2$ if $x_n<0$ 
and $R_-=0$ if $x_n=0$.  
Thus, we conclude that {\it if $n>6$, $x_1+x_n>0$, then $R_+=n-1, R_-\in\{0, n-2\}$. Moreover, $\mathfrak{q}_x=\mathfrak{q}_s$ where $s$ equals $(1, 0,\ldots, 0)$ or $(2,0,\ldots,0,-1)$.} \\

{\it Case 2:}  Suppose that $x_1+x_n<0$.  Then $R_-\ge n-1$. Since $n-1\ge R_+\ge R_-$ by hypothesis, it follows that $R_+=n-1=R_-$.
Now we proceed, using Lemma \ref{duality}, as in the previous case and conclude that if $n>6$, then there is no such $x$. 

{\it Case 3:} Suppose that $x_1+x_n=0.$  Set $a:=\#I_{[1]}, b:=\#J_{[1]}$.  We have $a,b\ge 1$. There are 
three subcases to consider:\\
{\it Subcase(i)}:  Let $a=1$.   {\it We claim that $b=1$ and that $R_\pm=n-2$ if $n\ge 7$.}  \\
Clearly $R_-\ge n-2$ as $x_i+x_n<0~\forall i>1$.  This implies that $R_+\ge n-2$.  
If $b>1$, then $x_i+x_n<0, x_j+x_{n-1}<0~\forall 2\le i\le n-1, 2\le j\le n-2$ and so 
$R_-\ge 2n-5$ implies that $n-1\ge 2n-5$ and 
so $n\le 4$, contrary to our hypothesis.  Hence we must have $b=1$.  
Let $R_+=n-1$ and $n\ge 7$.  We will arrive at a contradiction. 
Since $x_1+x_n=0, a=b=1, R_+\le n-1$, we must have $x_2+x_3\ge 0, x_2+x_j\le 0$ for $4\le j\le n$ and $x_i+x_j\le 0$ 
for $3\le i\le n$.  
Since $R_-\le n-1,$ we must have $x_{n-3}+x_{n-1}=0$ which implies that 
$x_2+x_4=0, x_2+x_5=0, x_3+x_4=0, x_3+x_5=0, x_4+x_5=0$. (For the last equality we have used $n-1>5$.) 
Therefore $x_2=x_3=x_4=x_5=0$.  In particular $x_2+x_3=0$ and so $R_+=n-2$, establishing our claim. \\ 

{\it Subcase (ii)}:  Let $b=1$. 
If $a>1$, again arguing as above we have $n-1\ge R_+\ge 2n-5$, which implies $n\le 6$.   So if $n\ge 7$ we are reduced to the 
previous subcase, namely $a=1$, and so conclude that $R_+=n-2=R_-$.\\

{\it Subcase (iii)}:  Suppose that $a, b>1$. 
Firstly, we have $a+b\le n$. Since $a, b\ge 2$, we have $x_i+x_j>0, 1\le i\le 2, ~i<j\le n-b$.  Therefore $n-1\ge R_+\ge n-1-b+n-2-b$. 
Hence $b\ge (n-2)/2$.  Write $n=2m+1$ or $n=2m$, according as $n$ is odd or even, so that $b\ge m$ or $m-1$. Similarly  $a\ge m$ 
or $m-1$, depending on the parity of $n$.

Suppose that $n=2m+1$.  Then 
Example \ref{minimalparabolicseg}(ii)  yields $n-1\ge R_+\ge a(n-a-b)+{m\choose 2}$.  So $2m\ge (m^2-m)/2$, which implies that $m\le 5$ and so $n\le 11$. 
Suppose that $n=11$.  If $a+b<11$, then $a=b=m=5$ and so we have the estimate $R_+\ge a(n-a-b)+{m\choose 2}=15.$ If $a+b=11$, 
then one of them, say $a$ equals $6$.  In this case $R_+\ge {a\choose 2}=15$.  Similarly, the possibility that 
$n=9$ is also eliminated.  

Suppose that $n=2m$.  Proceeding as in the case $n$ odd, we obtain that $a,b\ge m-1$, $2m-1\ge (m-1)(m-2)/2$ which implies $m\le 6$. When $n=10$, it is readily verified that $9\ge a(n-a-b)+{a\choose 2}$ has no solution when $a, b\ge 4$.  Similarly, 
when $n=12$, the inequality $11\ge a(n-a-b)+{a\choose 2} $ has no solution when $a,b\ge 5$. 
This establishes our claim.
\end{proof}

\begin{remark}{\em 
When $4\le n\le 8$, there are more possibilities for the $\theta$-stable parabolics with $R_+, R_-\le n-1$.  
The following list of sequences gives the complete list of the `exceptional' $\theta$-stable parabolic subalgebras 
$\mathfrak{q}
=\mathfrak{q}_x$  where $1\le R_-\le R_+\le n-1$,   
and the values of $R_+(\mathfrak{q}_x)$.   
In all these cases $R_+=R_-$.\\    
$\underline{n=8}$:  The only exceptional $\mathfrak{q}$ corresponds to $x=(1,1,1,1,-1,-1,-1,-1)$ where $R_+=6$.\\
$\underline{n=7}$:  The only exceptional $\mathfrak{q}$ corresponds to  $x=(1,1,1,0,-1,-1,-1)$, where $R_+=6$.\\ 
$\underline{n=6}$:  There two exceptions corresponding to $x=(1,1,1,-1,-1,-1)$ in which case $R_+=3$, 
and to $x=(2,2,1,-1,-2,-2)$ where $R_+=5$.\\
$\underline{n=5}$:  There are two exceptional cases 
corresponding to $x=(1,1,0,-1,-1)$ with $R_+=3; x=(2,1,0,-1,-2)$ with $R_+=4$.  \\
$\underline{n=4}$:  The only exceptional case corresponds to $x=(1,1,-1,-1)$ with $R_+=1$.
}
\end{remark}

\begin{remark} \label{equivalentq}
{\em 
(i)  Suppose that $\Phi(\mathfrak{q})$ is the set of weights of $\mathfrak{q}=\mathfrak{q}_x$.  Choose a positive system for 
$(\mathfrak{k},\mathfrak{t})$ with respect to which the roots in $\Phi_x^+\cap \Phi_\mathfrak{k}$ and the weights of $\mathfrak{u}\cap \mathfrak{k}$ are positive.  Such a positive system of $(\mathfrak{k},\mathfrak{t})$ is said to be compatible with $\Phi(\mathfrak{q})$.
In what follows, we talk of highest weight of a $K$-representation with respect to this positive system.   
Let $\mu=\mu(\mathfrak{q})$ be the sum of all the $\mathfrak{t}$-weights of $\mathfrak{u}\cap \mathfrak{p}$.  Then it 
is known that the irreducible $K$-representation $V_\mu$ with highest weight $\mu$ occurs in 
the $(\mathfrak{g},K)$-module $A_{\mathfrak{q},K}$ with multiplicity one.  
 Any other $K$-type occurring in $A_{\mathfrak{q},K}$ has highest weight equal to a sum 
 $\mu+\sum a_\beta \beta$ where $a_\beta\ge 0$ and 
$\beta$ varies over the weights of $\mathfrak{u}\cap \mathfrak{p}$.  Thus $V_\mu$ is the lowest $K$-type occurring in $A_\mathfrak{q}$.
Moreover, if $\mathfrak{q}':=\mathfrak{q}_{x'}$ and the positive system of $(\mathfrak{k},\mathfrak{t})$ is compatible with $\Phi(\mathfrak{q}')$, then $\mathcal{A}_{\mathfrak{q}}$ and $\mathcal{A}_{\mathfrak{q}'}$ are unitarily equivalent if and only if 
$\mathfrak{u}\cap \mathfrak{p}=\mathfrak{u}'\cap \mathfrak{p}$.   This is a particular case of a very general statement proved in \cite[Proposition 4.5]{sriba}. 

(ii)  In view of the fact that the $(\mathfrak{g},K)$-cohomology of $A_{\mathfrak{q},K}$ 
depends only on 
$R(\mathfrak{q})$ and the compact symmetric space $Y_\mathfrak{q}$ associated to $\mathfrak{q}$, 
it is possible that the Harish-Chandra modules associated to two non-isomorphic 
representations $\mathcal{A}_\mathfrak{q}, \mathcal{A}_{\mathfrak{q}'}$ can have isomorphic $(\mathfrak{g},K)$-cohomology.  
For example, if $\mathcal{A}_\mathfrak{q}$ is a discrete series representation we have $R(\mathfrak{q})=(1/2)\dim (G/K)$ and $H^*(\mathfrak{g},K;A_{\mathfrak{q},K})\cong \mathbb{C}$.  In this case, however, distinct such discrete series have 
distinct pairs $(R_+(\mathfrak{q}), R_-(\mathfrak{q}))$.   It is easy to construct examples of inequivalent 
representations whose $(\mathfrak{g},K)$-cohomologies have isomorphic Hodge types for {\it all} pairs $(p,q)$.  For 
example, this happens when $x=(2, 1, 1, -1,-1,-2)$, $x'=(2,2,1,-1,-2,-2)$.

(iii) Proposition \ref{minimalparabolics}, among many others, has been obtained by an entirely 
different approach, by Arghya Mondal \cite{mondal} using a combinatorial model referred to as the 
{\it decorated staircase diagram.} 
}
\end{remark}
\section{Proofs of Theorems \ref{main1} and \ref{main2}}\label{proofs}
We are now ready to prove the main results of the paper.  
If $\mathfrak{q}$ is a $\theta$-stable parabolic subalgera, we shall denote by $[\mathfrak{q}]$ the equivalence 
class of $\mathfrak{q}$ where $\mathfrak{q}\sim\mathfrak{q}'$ if the irreducible $G$-representations 
$(\mathcal{A}_\mathfrak{q},A_\mathfrak{q})$ and $(\mathcal{A}_{\mathfrak{q}'},A_{\mathfrak{q}'})$ are unitarily equivalent.  
When $\mathfrak{q}\sim\mathfrak{g},$ we have $\mathfrak{q}=\mathfrak{g}$ and the representation $\mathcal{A}_\mathfrak{q}$ is the trivial (one-dimensional) representation. 

\noindent
{\it Proofs of Theorem \ref{main1} and \ref{main2}:}  Let $u\in H^{p,p}(X_\Gamma;\mathbb{C})$.  Recall the Matsushima isomorphism 
$H^{p,p}(X_\Gamma)=\oplus_{[\mathfrak{q}]}m(\mathfrak{q},\Gamma)H^{p,p}(\mathfrak{g},K;A_{\mathfrak{q},K})$.  The summand 
corresponding to $\mathfrak{q}=\mathfrak{g}$ is the image of the Matsushima homomorphism $H^*(X_u;\mathbb{C})\to 
H^*(X_\Gamma;\mathbb{C})$.   
Write $u=\sum_{[\mathfrak{q}]} u_{[\mathfrak{q}]}$ where $u_{[\mathfrak{q}]}\in H^*(\mathfrak{g},K;A_{\mathfrak{q},K})$.  Then 
$u_{[\mathfrak{q}]}=0$ for all $\mathfrak{q}\ne \mathfrak{g}$  if and only if  
$u$ is in the image of the Matsushima homomorphism.  

Any sign involution or an involution of even type $\sigma$ (constructed in \S \ref{commutinginvolutions} ) defines an analytic cycle $C_\Lambda(\sigma)$ in the 
K\"{a}hler manifold
$X_\Lambda$ and hence a non-zero cohomology class $[C_\Lambda(\sigma)]\in H^{p,p} (X_\Lambda;\mathbb{C})$ where $p$ is the 
(complex) codimension of $C_\Lambda(\sigma)\subset X_\Lambda$.   See \S \ref{specialcycle}.  Since the cohomology class 
represented by a special cycle is not in the image of the Matsushima homomorphism, it does not belong to the component $H^*(\mathfrak{g},K;\mathbb{C})\subset H^*(X_\Lambda;\mathbb{C})$. 
It follows that $[C_\Lambda(\sigma)]_{[\mathfrak{q}]}\ne 0$ for some $\mathfrak{q}\ne \mathfrak{g}$.  For any such $\mathfrak{q}$,  
we must have $R_+(\mathfrak{q})=R_-(\mathfrak{q})\le p$.    

Taking $C:=C_\Lambda(\sigma)$ to be as in Example \ref{extraspecialcycle} we have $p=n-1$. It follows that  
$[C]_{[\mathfrak{q}_0]}\ne 0$ for some $\mathfrak{q}_0$ of type $(R_+,R_+)$ with $R_+\le n-1$.  
Hence $m(\mathfrak{q}_0,\Lambda)\ne 0$.   Since $Y_{\mathfrak{q}_0}$ is Hermitian symmetric by Lemma \ref{derivedl0}, we see that $H^{j,j}(Y_{\mathfrak{q}_0};\mathbb{C})\ne 0$ for $0\le j\le \dim_\mathbb{C}Y_{\mathfrak{q}_0}$.
Since $H^{r,r}(\mathfrak{g},K;A_{\mathfrak{q}_0,K})
\cong H^{r-R_+,r-R_+}(Y_{\mathfrak{q}_0};\mathbb{C})\ne 0$ for $R_+\le r\le \dim_\mathbb{C} X_\Lambda-R_+={n\choose 2}-R_+, $ 
we conclude that there are cohomology classes of Hodge type $(r,r)$  in 
$X_\Lambda$ which are not in the image of the Matsushima map for all $r$ such that $R_+\le r\le {n\choose 2}-R_+.$

When $n$ is odd we observed in \S \ref{specialcycle} that the dual cohomology class of a special cycle $[C(\sigma)]$ 
where $\sigma$ is of odd type is a non-zero cohomology class in the middle dimension ${n\choose 2}$.

 By Proposition \ref{minimalparabolics}, 
when $n>8$ and $1\le r\le n-1$, there is exactly one class of $\theta$-stable parabolic subalgebra $[\mathfrak{q}_0]$ of type $(r,r)$ 
corresponding to $r=n-2$.  So $[C_\Lambda(\sigma)]_{[\mathfrak{q}_0]}\ne 0$. 
In follows that $m(\mathfrak{q}_0,\Lambda)\ne 0$ for any lattice $\Lambda$ 
as in the statement of Theorem \ref{main1}.  From what has been established already, we see that there are cohomology classes 
of Hodge type $(r,r)$ for $n-2\le r\le {n\choose 2}-(n-2)$ which are not in the image of the 
Matsushima homomorphism $H^*(X_u;\mathbb{C})\to H^*(X_\Lambda;\mathbb{C})$.
This proves Theorem \ref{main1}.    \hfill $\Box$

\begin{remark}{\em 
Our method of proofs of the main results of this paper is applicable in a more general setting.  
Analogous results for the case $G=Sp(n,\mathbb{R})$ are under preparation. }
\end{remark}

\noindent
{\bf Acknowledgments}  The authors thank Prof. J. Oesterl\'e for his help with the proof of Lemma \ref{involutiontype}(ii).  Research of both authors were partially supported by the Department of Atomic Energy, 
Government of India, under a XII Plan Project.

\end{document}